\documentclass[12pt,a4paper]{article}

\usepackage{amssymb,amsmath,amsthm,xspace,epsfig, amsfonts,amscd}

\usepackage[left=1.5cm, right=1.5cm, top=1.5cm]{geometry}

\begin{document}

\newcommand{\N}{\mathbb{N}}
\newcommand{\R}{\mathbb{R}}
\newcommand{\Z}{\mathbb{Z}}
\newcommand{\Q}{\mathbb{Q}}
\newcommand{\C}{\mathbb{C}}
\newcommand{\PP}{\mathbb{P}}

\newcommand{\LL}{\Bbb L}
\newcommand{\OO}{\mathcal{O}}

\newcommand{\esp}{\vskip .3cm \noindent}
\mathchardef\flat="115B

\newcommand{\lev}{\text{\rm Lev}}

\def\ut#1{$\underline{\text{#1}}$}
\def\CC#1{${\cal C}^{#1}$}
\def\h#1{\hat #1}
\def\t#1{\tilde #1}
\def\wt#1{\widetilde{#1}}
\def\wh#1{\widehat{#1}}
\def\wb#1{\overline{#1}}

\def\restrict#1{\bigr|_{#1}}

\def\hu#1#2{\mathsf{U}_{fin}\bigl({#1},{#2}\bigr)}
\def\ch#1#2{\left(\begin{array}{c}#1 \\ #2 \end{array}\right)}

\newtheorem{lemma}{Lemma}[section]

\newtheorem{thm}[lemma]{Theorem}
\newtheorem*{thm*}{Theorem}
\newtheorem*{lemma?}{Lemma ??}

\newtheorem{defi}[lemma]{Definition}
\newtheorem{conj}[lemma]{Conjecture}
\newtheorem{cor}[lemma]{Corollary}
\newtheorem{prop}[lemma]{Proposition}
\newtheorem{claim}{Claim}
\newtheorem{prob}[lemma]{Problem}
\newtheorem{qu}[lemma]{Question}
\newtheorem{q}[lemma]{Questions}
\newtheorem*{rem}{Remark}
\newtheorem{examples}[lemma]{Examples}
\newtheorem{example}[lemma]{Example}

\title{Non-metrizable manifolds and contractibility}
\date{\today}
\author{Mathieu Baillif}
\maketitle

\abstract{\footnotesize 
   We investigate whether non-metrizable manifolds in various classes can be homotopy equivalent to a CW-complex (in short: heCWc),
   and in particular contractible.
   We show that 
   a non-metrizable manifold cannot be heCWc if it has one of the following
   properties: it contains a countably compact non-compact subspace; it contains a copy of an $\omega_1$-compact
   subset of an $\omega_1$-tree; it contains a non-Lindel\"of closed subspace 
   functionally narrow in it. (These results hold for more general spaces than
   just manifolds.)
   We also show that the positive part of the tangent bundle of $\LL_+$ is not heCWc (for any smoothing).
   These theorems follow from stabilization properties of real valued maps.
   On a more geometric side, we also show that the Pr\"ufer surface, which has been shown to be contractible
   long ago,
   has an open submanifold which is not heCWc.
   On the other end of the spectrum, we show 
   that there is a non-metrizable contractible Type I surface.
   }

\section{Introduction}

By {\em space} we mean `topological Hausdorff space', in particular `regular' and 'normal' imply Hausdorff.
By {\em manifold}, or {\em $n$-manifold} if we want to empasize the dimension,
is meant a connected space each of whose points has a neighborhood homeomorphic 
to the Euclidean
space $\R^{n}$ for some finite $n\ge 1$. 
A {\em surface} is a $2$-manifold. 
Every function is assumed continuous otherwise stated.
We use brackets $\langle\,,\,\rangle$ for ordered pairs reserving parenthesis $(\,,\,)$ for open intervals.
The letters $\alpha,\beta,\gamma$ are used exclusively for ordinals, and an ordinal is assumed to be the 
set of its predecessors.
Undefined terminology used in this introduction is explained either at its end or in later sections.

This article is concerned with homotopy questions for non-metrizable manifolds and similar spaces. 
This is 
a somewhat exotic subject for which there exists a small but growing literature;
let us cite the
book by D. Gauld \cite[Chapter 5.2]{GauldBook} 
and the articles by D.Gauld, S. Deo, A. Gabard and the author 
\cite{GabardWouuuh, meszigueshom, BaillifDeoGauld, DeoGauld:2007},
for instance. It is worth noting that the methods and results in these papers (and the present one)
belong mainly to general topology (with a small set theoretic flavor) and have nothing
in common with homotopy theory per se.

Acknowledging its relative unpopularity, let us explain why 
we think homotopy in non-metrizable manifolds is an interesting 
subject.\footnote{The true and honest reason is that we happen to like to think about these problems. 
It might unfortunately be seen as slightly insuficient.}
Milnor showed in \cite{MilnorCW} that a metrizable manifold 
is homotopy equivalent to a CW-complex (in short: {\em heCWc}), hence
metrizable manifolds that are contractible coincide with those having vanishing homotopy groups
by Whitehead Theorem (see, e.g., \cite[p. 346]{Hatcher}).
If one releases the metric assumption, this is not true anymore as
there are simple non-metrizable non-contractible manifolds with vanishing homotopy groups (e.g. the {\em longray} described below).
Of course, any space with a non-trivial homotopy group is non-contractible.
It was at first not clear to us whether contractible non-metrizable manifolds do exist at all.
We then noticed that examples were found a long time ago by Calabi and Rosenlicht:
the Pr\"ufer surface and some of its variants are contractible.
(The original source is \cite{CalabiRosenlicht}, see also \cite{GabardWouuuh}.) 

While manifolds can be widly different one from another, they share the local properties of the Euclidean space.
In particular, any manifold is locally compact and thus Tychonoff, locally (path) connected and locally second countable.
In addition, metric manifolds are normal 
(even perfectly normal and monotonically normal, see e.g. \cite[p. 23]{GauldBook} for definitions), 
second countable and thus Lindel\"of and paracompact. 
Lindel\"ofness and paracompactness are actually
two of the many properties equivalent to metrizability for manifolds, see \cite[Chapter 2]{GauldBook}.
In metric spaces, Lindel\"ofness is equivalent to being $\omega_1$-compact, and compactness to countable compactness.
(Recall that a space is {\em $\omega_1$-compact} or has {\em countable extent} iff its closed discrete subsets are at most countable.
A countably compact space is $\omega_1$-compact.)
Non-metrizable manifolds, on the contrary, can exhibit some but not all of these features, 
and the mere existence of non-metrizable manifolds with such-and-such properties sometimes depends on the axioms of set theory.
For example, perfectly normal manifolds do exist under {\bf CH} and do not under {\bf MA + $\neg$CH}
(due to M.E. Rudin and P. Zenor, see e.g. \cite[Chapter 6]{GauldBook}),
and monotonically normal non-metrizable
manifolds of dimension $>1$ do not exist (Z. Balogh and M.E. Rudin, see \cite[Cor. 2.3 (e)]{BaloghRudin:1992}).

Given this variety, we thought it might be interesting to see 
whether contractible non-metrizable manifolds exist in such-and-such class, or on the other hand 
if such-and-such topological property
prevents contractibility {\em regardless of the homotopy groups}. 
Some results are already available, for instance
the following theorem was proved by S. Deo and D. Gauld (and by the author for manifolds):

\begin{thm}[{\cite[Prop. 7.1]{meszigueshom}, \cite[Thm 3.4]{DeoGauld:2007}}]\label{thmdeogauld}
A locally metrizable space (in particular, a manifold) containing a copy of $\omega_1$ is non-contractible.
\end{thm}

Notice that the proper forcing axiom {\bf PFA} 
implies that a countably compact non-compact subset of a manifold contains a copy of $\omega_1$ 
(due to Balogh, see \cite{Balogh:1989}), an immediate corollary is that
under {\bf PFA} manifolds with such subsets are non-contractible. 
We show below that {\bf PFA} is actually not needed.

On the other side of the spectrum, the contractible
Pr\"ufer surface is neither normal nor $\omega_1$-compact.
(A version of it is realcompact, though, as shown by A. Mardani \cite{MardaniThesis}.)
It contains an open (Lindel\"of) Euclidean set with a non-metrizable (equivalently, non-Lindel\"of) closure.
In a terminology dating to P. Nyikos, it is in the 
broad class of
{\em Type II} manifolds (see just below for a precise definition). 
In a loose sense, what makes it non-metrizable
lies just at the boundary of some perfectly nice Euclidean open contractible set. We happen to be able to 
`push' this non-metrizable stuff inside this Euclidean open set,
yielding the contractibility.

The other broad class of non-metrizable manifolds are the {\em Type I} manifolds \cite[Def. 2.10]{Nyikos:1984}.
Denote as usual by $\wb{A}$ the closure of $A$ (in some topological space clear from the context).
A space $X$ is of Type I if and only if $X=\cup_{\alpha\in\omega_1}X_\alpha$, where $X_\alpha$ is open, $\wb{X_\alpha}$ Lindel\"of
for each $\alpha$, and $\wb{X_\alpha}\subset X_\beta$ whenever $\alpha<\beta<\omega_1$ (and of Type II otherwise).

In case of manifolds, each $X_\alpha$ is an open metrizable submanifold.
Loosely speaking again, these manifolds `grow slowly' 
instead of jumping at once into non-metrizability. In a Type I space, Lindel\"of subsets have Lindel\"of closure.

Let us give a summary of our results.
The first one (obtained with Peter Nyikos)
shows that Type I contractible manifolds exist.

\begin{thm}\label{thm:main-contractible}
   There is a contractible non-metrizable Type I surface.
\end{thm}

The mapping class group of such a manifold, on the other hand, can be quite large, see Section \ref{sec:mapping} 
for definitions
and details.
This surface is obtained by ``fattening'' an $\R$-special $\omega_1$-tree and is non-normal.
(Trees in this article are to be understood in a set theoretic sense. Definitions of the relevant notions
are recalled in Section \ref{sec:trees}.)
Our other results are on the negative side, we show that spaces in various classes are not contractible, 
actually not heCWc (which we recall is short for {\em homotopy equivalent to a CW-complex}).
We first obtained a generalization of Theorem \ref{thmdeogauld}.
(Theorems \ref{thm:main-cblycpct}--\ref{thm:main-dir} below are stated for manifolds in this
introduction, but we obtained results in larger classes.)

\begin{thm}\label{thm:main-cblycpct}
   A manifold containing a
   countably compact non-compact subspace is not heCWc.
\end{thm}
Since $\omega_1$-compactess is a weakening of countable compactness,
it seemed natural to see whether the same result holds in this class.
What we obtained is something quite weaker and involves subsets of trees.
The order on the tree is essential to our proof.

\begin{thm}\label{thm:main-w1cpct}
   A manifold containing a copy of an $\omega_1$-compact subset of an $\omega_1$-tree
   (in particular, a Suslin tree) is not heCWc.
\end{thm}

Recall that the closed longray $\LL_{\ge 0}$ is the product $\omega_1\times[0,1)$ 
with lexicographic order topology (not the product topology~!). When dealing with it,
we often see $\omega_1$ as a subset of $\LL_{\ge 0}$ by identifying $\alpha$ with $\langle \alpha,0\rangle$.
$\LL_+$ is $\LL_{\ge 0}$ with the $0$ point removed.

Our next result involves the concept of functionally narrow subspaces in Type I spaces,
which whose definition is recalled in Section \ref{sec:dir}.
Loosely speaking, a closed subspace in a Type I space $X$ is functionally narrow in $X$ if it cannot be
separated by a function $X\to\LL_{\ge 0}$ into two non-Lindel\"of pieces, one with a bounded image 
and the other with an unbounded image.
Functional narrowness is relative to the ambient space, however
a closed copy of $\omega_1$ in any Type I space is functionally narrow in it. 
There are spaces with discrete subspaces
functionally narrow in them as well, see \cite[Example 5.12]{mesziguesnarrow} (and the rest of this paper
for more on the subject).

\begin{thm}\label{thm:main-dir}
  A Type I (non-metrizable) manifold containing a closed non-Lindel\"of subspace functionally narrow in it is not heCWc.
\end{thm}

We also show that surfaces in a very interesting class studied extensively by P. Nyikos
in \cite{Nyikos:1992}
are non-contractible: the (positive part of the) tangent bundle of the longray $\LL_+$.
(We write {\em the} tangent bundle, but actually there are a lot of non-homeomorphic ones,
with distinct topological properties.)

\begin{thm}\label{thm:main-tangent}
   The positive part of the tangent bundle of $\LL_+$ for a given smoothing
   is not contractible (and thus not heCWc).
\end{thm}

(See the preprint \cite[Section 6]{mesziguessurf} for related results.)
The reason for the non-homotopy equivalence in Theorems \ref{thm:main-cblycpct}--\ref{thm:main-tangent} relies on
stabilization properties of real valued maps on relevant
subspaces.
This is basically the same for Theorem \ref{thmdeogauld}, whose proof is based on the fact
that any function $\omega_1\to\R$ is eventually constant.
However, 
a non-metrizable manifold can be non-contractible for 
more `geometrical' reasons: there is a submanifold (with trivial homotopy groups)
of the contractible Pr\"ufer surface which is non-contractible.
To understand better the statement, let us briefly recall that (the separable version of) 
Pr\"ufer surface
is obtained by adjoining to the open Euclidean 
plane uncountably many boundary components $\R_c$, each homeomorphic to $\R$.
(This version of Pr\"ufer surface is actually a surface {\em with boundary}, see 
the end of this introduction.)
A precise definition is recalled in Section \ref{sec:prufer}.

\begin{thm}\label{thm:main-subprufer}
   The Pr\"ufer surface minus the $0$-points in each $\R_c$ is a non-contractible
   (and hence not heCWc) non-metrizable submanifold of a contractible manifold.
\end{thm}

This last example gives us the impression that being heCWc
is a somewhat rare occurence in the world of non-metrizable manifolds: some kind of perfect
alignment of uncountably many pieces seems to be necessary.
We are thus inclined to believe that the problems below have negative answers, but
were not able to find any evidence corroborating our feeling.

\begin{prob}\label{prob:w1cpct}
   Is there a (consistent) example of a contractible $\omega_1$-compact manifold~?
   More generally, is there 
   a contractible locally metrizable space containing a non-metrizable $\omega_1$-compact subspace~? 
\end{prob}

Theorem \ref{thm:main-w1cpct} is far from providing an answer to this problem.
For instance, with the help of $\diamondsuit$ one can build an $\omega_1$-compact
$2$-manifold which mimicks any given $\omega_1$-tree, in particular an $\R$-special one, see
\cite{Greenwood:2002, GreenwoodNyikos:2005}. 
(Recall that the surface
in Theorem \ref{thm:main-contractible} is build with an $\R$-special tree.)
It is not clear to us whether these manifolds are contractible. 
Under {\bf PFA}, a contractible $\omega_1$-compact manifold cannot be of Type I (see Theorem \ref{thm:PFAw1cpct} below),
but this does not settle  the general case.
Another problem, for which we have no information so far:

\begin{prob}\label{prob:normal}
   Is there a normal contractible non-metrizable manifold~? A perfectly normal one, consistently~? 
\end{prob}

Notice that by Jones Lemma (see e.g. \cite[1.7.12 (c) p.60]{Engelking}), 
a separable normal manifold is $\omega_1$-compact if $2^\omega < 2^{\omega_1}$,
hence Problems \ref{prob:w1cpct}--\ref{prob:normal} might be related.

As said above, Theorems \ref{thmdeogauld} and \ref{thm:main-cblycpct}--\ref{thm:main-dir} 
hold for more general spaces than just manifolds.
It might be instructing to note that some amount of non-compactness, first countability and local Lindel\"ofness  
has to be involved to ensure non-contractibility,
as the next simple examples show.

\begin{example}\label{ex:cones}
   There are contractible non-metrizable spaces in the following classes:\\
   (i) compact and first countable,\\
   (ii) compact and containing a (non-closed) copy of $\omega_1$,\\
   (iii) first countable, countably compact and containing a closed copy of $\omega_1$.
\end{example}
\proof[Details]
For (i), Helly's space \cite[Ex. 107]{CEIT} is compact, first countable and contractible. 
For (ii)--(iii), given a space $X$ we 
let $CX=X\times[0,1]/\langle x,1\rangle\sim\langle y,1\rangle$ be the cone over $X$.
$CX$ is contractible for any space $X$. For (iii) take $C\omega_1$, which is first countable
(a neighborhood basis of the vertex of the cone is given by the vertex union 
$\omega_1\times (1-1/n,1)$ for $n\in\omega$, as well known).
Notice that $C\omega_1$ is not locally Lindel\"of. For (ii)
take $C(\omega_1+1)$. It is not first countable (or even countably tight, for that matter)
due to extra point in $\omega_1 + 1$. 
\endproof

\vskip .3cm
This paper is organized in somewhat independant sections.
In Sections \ref{sec:cblycpct} \& \ref{sec:dir}, we prove Theorems \ref{thm:main-cblycpct} and \ref{thm:main-dir}
involving countable compactness and functionally narrow subspaces; Section \ref{sec:tangent} deals with Theorem \ref{thm:main-tangent}
(about tangent bundles of $\LL_+$)
and Section \ref{sec:prufer} with Theorem \ref{thm:main-subprufer} (about subspaces of Pr\"ufer surface).
Then, Section \ref{sec:trees} is a bit different as it studies $\omega_1$-trees and the homotopic 
properties of their road spaces,
generalizing Theorem \ref{thmdeogauld} in various ways. This enables us in particular to prove Theorem \ref{thm:main-w1cpct}.
We also prove that the road space of an $\R$-special tree is contractible, which gives half the proof of 
Theorem \ref{thm:main-contractible}. The other half was given to us by Peter Nyikos who found a way
to ``fatten'' this road space to obtain a surface homotopy equivalent to it. This part is given in Section \ref{sec:surfaces}. 
This
Section contains further results on the mapping class group of this contractible surface and the Pr\"ufer surface.
(We note in passing that the results of Section \ref{sec:trees} appeared in the preprint \cite{meszigues-contract-trees}.)

\vskip .3cm
Let us finish this introduction with some definitions for further reference.
By {\em cover} of a space, we always mean a cover by open sets.
A {\em manifold with boundary} $M$ is a connected space each of whose points has a neighborhood homeomorphic 
to $\R^{n-1}\times \R_{\ge 0}$ for some finite $n>0$.\footnote{Whether a manifold possesses a 
boundary or not is mostly irrelevant for the
content of this paper, but we sometimes allow boundaries to simplify some definitions.}
Its {\em boundary} $\partial M$ 
contains the points which do not have a neighborhood homeomorphic to $\R^n$, it is an $n-1$-manifold (without
boundary) and a closed subset of $M$.
The term {\em club} is a shorthand for {\em closed and unbouded}.
We use very often the notation $h_t(x)$ for $h(x,t)$ and $h_t$ for $h(\cdot,t)$
when $h$ is a function with domain $X\times[0,1]$.
An {\em homotopy} between two maps $f,g:X\to Y$ is a function $h:X\times[0,1]\to Y$ such that $h(\cdot,0) = f$
and $h(\cdot,1) = g$ (in short: $h_0=f$, $h_1=g$). 
Two spaces $X,Y$ are {\em homotopy equivalent} if there are $f:X\to Y$ and $g:Y\to X$ such that
$f\circ g$ and $g\circ f$ are homotopic to $id_Y$ and $id_X$, respectively.
A {\em contraction} is an homotopy between the identity map and a constant map, and a space is {\em contractible}
if it admits a contraction (or equivalently if it is homotopy equivalent to a one point space).

We recall briefly the definition of a CW-complex.
(Caution: to avoid confusion with the notation $X^n$ usually reserved for the $n$-th power
of the space $X$ in point set topology, we will denote the $n$-skeleton
by $X^{(n)}$.)
A CW-complex is a space $X$ constructed by starting with a discrete space $X^{(0)}$ of $0$-cells 
(which are points).
We then inductively form the $n$-skeleton $X^{(n)}$ by attaching $n$-cells (which are copies of
the open $n$-disk) $e^n_i$ ($i$ in some index set)
via continuous maps
$\varphi_i:S^{n-1}\to X^{(n-1)}$; that is, we identify $x\sim \varphi_i(x)$ for $x\in\partial D^n_i$, where $D^n_i$ is a copy 
of the closed $n$-disk and $\partial D^n_i \simeq S^{n-1}$ its boundary. We endow $X^{(n)}$ with the quotient topology.
The $n$-cell $e^n_i$ is then the homeomorphic image of $D^n_i  -  \partial D^n_i$.
$X=\cup_{n\in\omega}X^{(n)}$ is given the weak topology: $A\subset X$ is open (or closed) 
if and only if $A\cap X^{(n)}$ is open (or closed) for all $n$.
In particular, $X^{(n)}$ is closed in $X$.
Then, $e^n_i$ 
is the image
under a characteristic map $\Phi_i:D^n_i\to X$ (the composition of $\varphi_i$ and the inclusion)
of $D^n_i -  \partial D^n_i$. 
For more details, see for instance Appendix A in \cite{Hatcher}, from which we borrowed the terminology.
Any CW-complex is hereditarily paracompact (see for instance \cite[Thm 1.3.5 \& Ex. 1 p. 33]{FritschPiccinini}) and 
hence hereditarily collectionwise normal (see \cite[Theorem 5.1.18 p. 305]{Engelking} for the definition and proof).

%%%%%%%%%%%%%%%%%%%%%%%%%%%%%%%%%%%%%%%%%%%%%%%%%%%%%%%%

\section{Countably compact non-metrizable manifolds are non-contractible}\label{sec:cblycpct}

A space is {\em C-closed} if any countably compact subspace is closed. 
The class of C-closed spaces
contains in particular the hereditarily paracompact
spaces, the sequential spaces and the regular spaces
with $G_\delta$ points, see e.g. \cite{IsmailNyikos}. In particular, manifolds and CW-complexes are C-closed.
A cover of a space is a {\em chain cover} if it is linearly ordered by the inclusion relation.
We may (and do) assume that chain covers are indexed by regular cardinals and that a member with a 
strictly smaller index than another member is strictly contained in it.
We will prove the following.

\begin{thm}\label{thm:cblycpct}
   Let $X$ be countably compact and $Y$ be C-closed such that each point in $Y$ has a Lindel\"of 
   (non-necessarily open) neighborhood.
   Let $h:X\times[0,1]\to Y$ be such that $h_{t_0}(X)$ is compact for some $t_0\in[0,1]$.
   Then $h_t(X)$ is compact for each $t\in[0,1]$.
\end{thm}

Recall that compactness is equivalent to linear compactness, that is, 
a space $S$ is compact if and only if any chain cover of it contains $S$ as a member.
(To our knowledge, this dates back to Alexandroff and Urysohn \cite{AlexandroffUrysohn}.)
Hence, if $S$ is a non-compact space, there is a chain cover $\mathcal{U}=\{U_\alpha\,:\,\alpha\in\kappa\}$ of $S$, with
$\kappa$ a regular cardinal, such that $S\not\subset U_\alpha$ for each $\alpha\in\kappa$.
Given such a chain cover $\mathcal{U}$, we say that $E\subset S$ is {\em $\mathcal{U}$-unbounded} iff
$E\not\subset U_\alpha$ for each $\alpha\in\kappa$.

\begin{lemma}\label{lemma:cUunbd}
   Let $X$ be countably compact and non-compact, and let $\mathcal{U}=\{U_\alpha\,:\,\alpha\in\kappa\}$ be a chain cover 
   such that $X\not\subset U_\alpha$ for each $\alpha$. 
   Then the following hold.
   \\
   (a) Given $f:X\to Y$ such that $f(X)$ is compact and $Y$ is C-closed,  
       there is $c\in Y$ such that $f^{-1}(\{c\})$ is $\mathcal{U}$-unbounded.\\
   (b) Let $E_n\subset X$ ($n\in\omega$) 
       be closed and $\mathcal{U}$-unbounded such that $E_{n+1}\subset E_n$. Then $E=\cap_{n\in\omega}E_n$ is closed 
       and $\mathcal{U}$-unbounded.
\end{lemma}
\proof
   \ \\
   (a) For each $\alpha\in\kappa$ choose $x_\alpha\in X-U_\alpha$ and set $E_\alpha = \{x_\beta\,:\,\beta\ge\alpha\}$.
   Then $\wb{E_\alpha} \cap U_\alpha = \varnothing$ and each $E_\alpha$ is $\mathcal{U}$-unbounded.
   Since $\wb{E_\alpha}$ is countably compact, $f(X)$ compact and $Y$ C-closed, then $f(\wb{E_\alpha})\subset f(X)$ 
   is compact. Any finite intersection of $E_\alpha$ is nonempty, hence
   by compactness $\cap_{\alpha\in\kappa} f(\wb{E_\alpha})\not=\varnothing$.
   Take $c$ is this intersection, then $f^{-1}(\{c\})\cap \wb{E_\alpha}\not=\varnothing$ for each $\alpha$, which shows that
   $f^{-1}(\{c\})$ is $\mathcal{U}$-unbounded.\\
   (b)
   For each $\alpha$, $\cap_{n\in\omega}(E_n-U_\alpha)\not=\varnothing$ by countable compactness.
\endproof

\proof[Proof of Theorem \ref{thm:cblycpct}]
   If $X$ is compact, there is nothing to do. We assume that $X$ is non-compact.
   Set $A = \{t\in[0,1]\,:\,h_t(X)\text{ is compact}\}$. We will show that $A$ is open and closed in $[0,1]$. Since
   $t_0\in A$, it is non-empty and thus equal to $[0,1]$.
   \\
   {\em $A$ is open.}
   Let $t\in A$. Assume that there is no open interval around $t$ in $A$, hence
   there is a sequence $t_n$ converging to $t$ such that $h_{t_n}(X)$ is non-compact for each $n\in\omega$.
   By local Lindel\"ofness of $Y$ and compactness of $h_t(X)$, there
   is a Lindel\"of $V\subset Y$ whose interior contains $h_t(X)$. Since $Y$ is C-closed and $h_{t_n}(X)$ is countably compact,
   the latter is closed for each $n$.
   Hence $h_{t_n}(X)\not\subset V$ (otherwise it would be compact), 
   choose $x_n\in X$ such that $h_{t_n}(x_n)\not\in V$.
   Take an accumulation point $x$ of the $x_n$'s, then $h_t(x)$
   is not in the interior of $V$, which contains $h_t(X)$, a contradiction. 
   \\
   {\em $A$ is closed.}
   Let $t_n$ be a sequence converging to $t$ and assume that $h_{t_n}(X)$ is compact for each $n$.
   We show that $h_t(X)$ is compact as well.
   Assume it is not the case, hence there is a chain cover 
   $\mathcal{V}=\{V_\alpha\,:\,\alpha\in\kappa\}$ of $h_t(X)$ such that $h_t(X)\not\subset V_\alpha$ for each $\alpha$.
   Let $\mathcal{U}$ be the chain cover of $X$ given by the inverse images of $V_\alpha$ under $h_t$.
   Set $E_0=X$. We assume by induction that $E_n$ is closed and $\mathcal{U}$-unbounded,
   in particular $E_n$ is countably compact non-compact.
   Since $h_{t_n}(E_n)$ is compact, by Lemma \ref{lemma:cUunbd} (a) there is $c_n\in Y$ such that
   $E_{n+1} = h_{t_n}^{-1}(\{c_n\})\cap E_n$ is $\mathcal{U}$-unbounded.
   By Lemma \ref{lemma:cUunbd} (b), $E=\cap_{n\in\omega}E_n$ is $\mathcal{U}$-unbounded.
   But $h_{t_n}$ is constant on $E$ for each $n$, hence by continuity $h_t$ is also constant on $E$, taking value $c\in f(X)$.
   Let $\alpha$ be such that $c\in V_\alpha$, then $E\subset U_\alpha=h_t^{-1}(V_\alpha)$, a contradiction.
\endproof
Notice that our only use of local Lindel\"ofness is to ensure that given any compact subspace of $X$,
there is a (non-necessarily open) neighborhood $V$
of it such that any countably compact subset of $V$ is compact. Notice also
that the proof goes through if one replaces $[0,1]$ by any connected sequential space.
Theorem \ref{thm:main-cblycpct} is an immediate consequence of Theorem \ref{thm:cblycpct}.

\proof[{Proof of Theorem \ref{thm:main-cblycpct}}]
   Manifolds are C-closed and each point has a Lindel\"of (open) neighborhood. 
   Let thus $M$ be a manifold containing a countably compact non-compact subspace $N$.
   Let $W$ be a CW-complex, and $f:M\to W$, $g:W\to M$ be given.
   Since $f(N)$ is countably compact and CW-complexes are paracompact and C-closed, $f(N)$ is compact,
   hence so is $g\circ f(N)$.
   It follows that $g\circ f$ cannot be homotopic to $id_M$ by Theorem \ref{thm:cblycpct}. 
\endproof

We close this section with a consequence of Theorem \ref{thm:cblycpct} under 
the locally compact trichotomy axiom {\bf LCT} due to Eisworth and Nyikos \cite{EisworthNyikos}.
Recall that a space is {\em $\omega$-bounded} iff any countable set has compact closure.
{\bf LCT} states that a locally compact space is either a countable union
of $\omega$-bounded spaces, or contains an uncountable closed discrete set,
or has a countable subset with non-Lindel\"of closure.
This axiom
holds for instance under {\bf PFA} or in models of the form {\bf PFA(S)[S]}, see
\cite{NyikosZdomskyy}.

\begin{thm}[{\bf LCT}]
   \label{thm:PFAw1cpct}
   Let $X$ be locally compact, C-closed, $\omega_1$-compact and non-Lindel\"of.
   If $X$ is heCWc, then $X$ does not contain a closed Type I subspace.
\end{thm}
\proof 
   Recall that we assume our Type I spaces to be non-Lindel\"of.
   It is well known that in a locally compact space, each point has an open Lindel\"of neighborhood.
   If $Y\subset X$ is closed and of Type I, the closure (in $Y$, hence in $X$) of any Lindel\"of subset of $Y$ is Lindel\"of.
   Hence by 
   {\bf LCT} $Y$
   is a countable union of countably compact subsets $Y_n$, which are thus closed in $X$.
   Since one of the $Y_n$ must be non-compact, 
   Theorem \ref{thm:cblycpct} shows that $X$ is not heCWc.
\endproof

%%%%%%%%%%%%%%%%%%%%%%%%%%%%%%%%%%%%%%%%%%%%%%%%%%%%%%%%%%%%%%%%%%%%%%%%%%%
%%%%%%%%%%%%%%%%%%%%%%%%%%%%%%%%%%%%%%%%%%%%%%%%%%%%%%%%%%%%%%%%%%%%%%%%%%%

\section{A Type I space containing a closed non-Lindel\"of subspace functionally narrow in it is non-contractible}\label{sec:dir}

Let $X$ be of Type I, hence $X=\cup_{\alpha\in\omega_1}X_\alpha$ with
$X_\alpha$ open, $\wb{X_\alpha}$ Lindel\"of and $\wb{X_\alpha}\subset X_\beta$ whenever $\alpha<\beta$.
The cover by the $X_\alpha$'s is said {\em canonical} if $X_\alpha = \cup_{\beta<\alpha}X_\beta$ for each limit $\alpha$.
Each chain cover indexed by ordinals $<\omega_1$ can be made canonical by adding members whenever needed.
(Note: our terminology differs slightly from that of Nyikos in \cite{Nyikos:1984} who uses the terms {\em canonical sequence} instead
of canonical cover.)
Two canonical covers agree on a club set of $\alpha$'s, as easily seen.
A subset of $X$ is {\em bounded} if it is contained in some $X_\alpha$ and {\em unbounded} otherwise.
A map between Type I spaces is {\em bounded} if its range is bounded and unbounded otherwise.
A {\em slicer} of (the canonical cover of) $X$ is a map $\sigma:X\to\LL_{\ge 0}$
such that $\sigma(\wb{X_{\alpha+1}}-X_\alpha)\subset[\alpha,\alpha+1]\subset\LL_{\ge 0}$.
If $X$ is regular, then $X$ is Tychonoff, each $\wb{X_\alpha}$ is normal and there is a slicer of $X$ 
(see \cite[Lemmas 2.2--2.3]{mesziguesnarrow}).
The following is proved in \cite[Theorem 3.8]{mesziguesnarrow}.
\begin{lemma}\ \\
  \label{lemmafrog}
   Let $X=\cup_{\alpha<\omega_1}X_\alpha$ be Type I and $D\subset X$ be club. Then the following are equivalent.\\
   (a) For each $f:X\to\LL_{\ge 0}$ with $f(D)$ unbounded, 
  the set $S\subset\omega_1$ of $\alpha$ such that $f(D -  X_\alpha)\subset[\alpha,\omega_1)$ is club,\\
   (b) For each $f:X\to\LL_{\ge 0}$ with $f(D)$ unbounded, 
  $f^{-1}(\{0\})\cap D$ is bounded in $X$.
\end{lemma}
A closed subspace $D$ of $X$ is {\em functionally narrow in $X$} if it satisfies (a) and (b) in the above lemma.
Whether a Type I space $X$ contains a non-Lindel\"of
subset functionally narrow in it can be a tricky question: any closed unbounded copy of $\omega_1$
is functionally narrow in $X$, but there are spaces containing discrete 
subspaces functionally narrow in them as well, and an Aronszajn tree (see the definition
in Section \ref{sec:trees}) does not
contain any uncountable subspace functionally narrow in it. What we prove is the following.

\begin{thm}\label{thm:dir}
   A regular Type I space containing a closed non-Lindel\"of subset functionally narrow in it is not heCWc.
\end{thm}

For the proof, we need a fact about CW-complexes.

\begin{lemma}\label{YuncountableinCW}
  Let $W$ be a CW-complex, and $Y=\{y_\alpha\, : \, \alpha\in\omega_1\}$ be an uncountable subset of $W$ 
  (all $y_\alpha$ are distinct). Then, either there is a cell $e^n_i$ such that
  $\wb{e^n_i}\cap Y$ is uncountable, or there is an uncountable $E\subset\omega_1$ and for each $\alpha\in E$,
  an open $V_\alpha\ni y_\alpha$ such that $\{V_\alpha\,:\,\alpha\in E\}$ is a discrete family.
\end{lemma}
Recall that a family $\mathcal{F}$ of subsets of $W$ is {\em discrete} iff any point in $W$ has a neighborhood which intersects at most 
one member of $\mathcal{F}$. This implies in particular that for any $\mathcal{G}\subset\mathcal{F}$,
$\cup\{\wb{F}\,:\,F\in\mathcal{G}\}$  is closed.

\proof
  We can assume that $Y\subset W^{(n)}$ for the smallest $n$ such that $W^{(n)}\cap Y$ is uncountable.
  Then, $Y\cap W^{(n-1)}$ is at most countable, so we can assume that $Y\subset W^{(n)} -  W^{(n-1)}$. So
  each $y_\alpha$ lies in an $n$-cell.
  If $e^n_i\cap Y$ is at most countable for all $n$-cells $e^n_i$,   
  take an uncountable subset $E\subset\omega_1$ such that each $e^n_i$ intersects $Y_0=\{y_\alpha\,:\,\alpha\in E\}$
  in at most one point. 
  It thus follows that $Y_0$ is closed discrete. 
  Since any CW-complex is collectionwise normal, it implies in particular
  that $Y_0$ can be expanded into a discrete collection of open sets, which gives us
  the required $V_\alpha$. 
\endproof

\proof[Proof of Theorem \ref{thm:dir}.]
  Let $X=\cup_{\alpha\in\omega_1}X_\alpha$ be a Type I regular space and $D\subset X$ be 
  closed, non-Lindel\"of and functionally narrow in $X$.
  Let $\sigma:X\to\LL_{\ge 0}$ be a slicer
  and $h:X\times[0,1]\to X$ be such that $h_1=id$.
  We show that $h_t(D)$ must be unbounded for each $t$.
  Set 
  $$
    \tau=\inf\{t\,:\,h_{s}(D)\subset X \text{ is unbounded }\forall s> t\}.
  $$ 
  We first show that $h_\tau$ must be unbounded on $D$. This is obvious if $\tau = 1$. If $\tau<1$,
  let $\{t_n\}_{n\in\omega}$ be a dense subset of $(\tau,1]$. By Lemma
  \ref{lemmafrog} (a) applied to $\sigma\circ h_{t_n}$, there is a club $S_n\subset\omega_1$
  such that $h_{t_n}(D -  X_\alpha)\subset X -  X_\alpha$ for all $\alpha\in S_n$.
  By continuity, this holds also for $h_\tau$ and $\alpha\in S=\cap_{n\in\omega}S_n$. Since
  $S$ is also club, $h_\tau(D)$ is unbounded. 
  If $\tau >0$,
  take a sequence $t_n\nearrow\tau$ with $h_{t_n}(D)$ bounded. By continuity, $h_\tau(D)$ must also be bounded.
  This shows that $\tau = 0$, hence $h_t$ is unbounded for each $t$.
  \\
  Now, for each $\alpha\in\omega_1$ choose $x_\alpha \in D -  X_\alpha$.
  Let $W$ be a CW-complex, and 
  $f:X\to W$, $g:W\to X$ be homotopy equivalences. Set $Y=\{y_\alpha=f(x_\alpha)\,:\,\alpha\in\omega_1\}$. If 
  $Y$ is countable or there is an uncountable $Y'\subset Y$  
  contained in a countable (hence Lindel\"of) subcomplex,
  then an unbounded subset of $D$ is sent by $g\circ f$ to a bounded subset of $X$.
  Since $D$ is functionally narrow in $X$,
  $\sigma\circ g\circ f(D)$ must be bounded by Lemma \ref{lemmafrog} (recall that $\sigma$ is a slicer), 
  thus so is $g\circ f(D)$. It follows that 
  $g\circ f$ cannot be homotopic to the identity. 
  \\
  Thus, we can assume that no uncountable subset of $Y$ is contained in a countable subcomplex of $W$.
  Let $E\subset\omega_1$ and the $V_\alpha\ni y_\alpha$ be given by Lemma \ref{YuncountableinCW}, and partition $E$
  into two disjoint uncountable subsets $E_0, E_1$. 
  Letting $U_\alpha=f^{-1}(V_\alpha)$, by discreteness the $\wb{U_\alpha}$ are $2$-by-$2$ disjoint and
  $\wb{\cup_{\alpha\in E_0}U_\alpha}=\cup_{\alpha\in E_0}\wb{U_\alpha}$. 
  We want to define a continuous
  $\chi:X\to \LL_{\ge 0}$ 
  such that $\chi(x_\alpha)=\alpha$ when $\alpha\in E_0$ and with value $0$ on the complement of $\wb{\cup_{\alpha\in E_0}U_\alpha}$.  
  This can be done since we assumed $X$ to be regular and hence (since of Type I), Tychonoff: 
  set $\chi$ to be $0$ on $X -  \wb{\cup_{\alpha\in E_0}U_\alpha}$,
  and inside $\wb{U_\alpha}$ take
  any function into $[0,\alpha]\subset\LL_{\ge 0}$ 
  that maps $x_\alpha$ to $\alpha$ and $\wb{U_\alpha} -  U_\alpha$ to $0$.
  The inverse image of an interval $(a,b)\subset\LL_{\ge 0}$ 
  (resp. $[0,b)\subset\LL_{\ge 0}$) is then an union of open subsets of the $U_\alpha$ (resp.
  together with $X - \cup_{\alpha\in E_0}\wb{U_\alpha}$), and is then open.
  By definition, $\chi(D)$ is unbounded (because of $E_0$), and $\chi^{-1}(\{0\})\cap D$ is unbounded as well
  (because of $E_1$). By Lemma \ref{lemmafrog}, $D$ is not functionally narrow in $X$, a contradiction.
  \endproof

%%%%%%%%%%%%%%%%%%%%%%%%%%%%%%%%%%%%%%%%%%%%%%%%%%%%%%%%%%%%%%%%%%

\section{Tangent bundles of $\LL_+$ are non-contractible}\label{sec:tangent}

We will use the terminology for Type I spaces defined at the beginning of the previous section.
The following theorem shows that a topologically very interesting
class of surfaces with trivial homotopy groups does not contain contractible ones.
Given a smoothing of the longray $\LL_+$, write $TL^+$ for the resulting tangent bundle.
P. Nyikos has extensively studied these spaces in \cite{Nyikos:1992}
and showed that they have very different topological properties depending on the smoothing.
We refer to this paper for precise definitions. Since we will just use properties shared by
all tangent bundles of $\LL_+$, we only need to recall
that there is a bundle projection $\pi: TL^+\to\LL_+$ such that $\{\pi^{-1}((0,\alpha))\,:\,\alpha\in\omega_1\}$
is a canonical cover for $TL^+$ 
and $\pi^{-1}((0,\alpha))$ is homeomorphic to $\R^2$ for each $\alpha$, hence $TL^+$ is a Type I surface.
(Here and in the remainder of this section, we view $\omega_1$ as a subset of $\LL_+$ whenever convenient.)
$TL^+$ contains a copy of $\LL_+$ (the $0$-section). Removing it
disconnects the surface in two homeomorphic subsurfaces $T^+,T^-$, none of which contain a copy of $\LL_+$.
Depending on the smoothing, they may contain a copy of $\omega_1$.
It is important to note that $\pi^{-1}((0,\alpha))\cap T^+$ is also homeomorphic to $\R^2$, thus
we slightly abuse the notation by denoting the restriction of $\pi$ to $T^+$ also by $\pi$.
We say that a function $f$ with domain $T^+$ is {\em almost constant}
if there is some club $C$ of $\omega_1$
such that $f$ is constant on $\pi^{-1}(C)\subset T^+$.

\begin{thm} 
   \label{tangentnotcontractible}
   If $h:T^+\times[0,1]\to\LL_+$  is given such that $h_{t_0}(T^+)$ is Lindel\"of for some $t_0\in[0,1]$,
   then $h_t$ is almost constant for each $t\in[0,1]$.
\end{thm}

An immediate corollary is:
\begin{cor}\label{cor:tangent}
    Neither $TL^+$ nor $T^+$ are heCWc. 
\end{cor}
\proof
The homotopy groups of $T^+$ are all trivial since any compact subset is contained in $\pi^{-1}((0,\alpha))$
(which is contractible)
for some $\alpha$.
If $T^+$ was heCWc it would then be contractible,
which would yield an homotopy between $\pi:T^+\to\LL_+$ and a constant map.
Theorem \ref{tangentnotcontractible} impedes it, proving the result for $T^+$.
The result for $TL^+$ follows from 
Theorem \ref{thmdeogauld} since it contains a copy of $\omega_1$.
\endproof

Following \cite{Nyikos:1992}, call a subset of $T^+$ {\em large}
if it intersects all the fibers above some stationary subset of $\omega_1$. 
We are going to use the following two results of P. Nyikos.
\begin{lemma}[{\cite[Thm 4.10 \& Cor. 4.15]{Nyikos:1992}}]
   \label{lemma:Nyikossmooth}
   \ \\
   (i) If $U$ is a large open set of $T^+$, then $\wb{U}$ contains $\pi^{-1}(C)$ for a club subset $C$ of $\omega_1$.\\
   (ii) If $U\subset T^+$ is open and contains $\pi^{-1}(C)$
        for some club $C$ of $\omega_1$, then given 
        $f:U\to\R$, there is a club $D$ of $\omega_1$ such that $f$ is constant on $\pi^{-1}(D)\cap U$. 
\end{lemma}
Actually, Nyikos proves part (ii) only for $U=T^+$, but his proof works as well in our case.

\proof[Proof of Theorem \ref{tangentnotcontractible}]
Let $h$ be as in the statement.
Write $T^+ = \cup_{\alpha\in\omega_1} T^+_\alpha$ with $T^+_\alpha=\pi^{-1}((0,\alpha))$ homeomorphic to $\R^2$.
Set 
$$ A = \{t\in[0,1]\,:\, h_t\text{ is almost constant}\}.$$
We show that $A$ is open and closed in $[0,1]$. 
Since $h_{t_0}(T^+)$ is Lindel\"of, it is contained in some $(0,\alpha)\subset\LL_+$ which is homeomorphic to $\R$,
hence $t_0\in A$ by Lemma \ref{lemma:Nyikossmooth} (ii). It follows that $A=[0,1]$.\\
{\em $A$ is closed.} Let $t_n\in A$ ($n\in\omega$) be a sequence converging to $t\in[0,1]$ and let $C_n\subset\omega_1$ 
be club
such that $h_{t_n}$ is constant on $\pi^{-1}(C_n)$. Then $C = \cap_{n\in\omega}C_n$ is club and 
by continuity $h_t$ is constant on $\pi^{-1}(C)$.\\
{\em $A$ is open.}
Let $t\in A$ and suppose in the way of contradiction that there is a sequence $t_n\not\in A$ ($n\in\omega$) converging to $t$.
Let $\gamma$ be such that $h_t(\pi^{-1}(C))$ is contained in $(0,\gamma)\subset\LL_+$ for some club $C$.
If $\pi\bigl(h_{t_n}^{-1}((0,\gamma))\bigr)$ is stationary, 
by Lemma \ref{lemma:Nyikossmooth} (i) $h_{t_n}^{-1}((0,\gamma+1))$ contains $\pi^{-1}(D_n)$ for a club $D_n$,
hence by (ii) $t_n\in A$. There is thus a club
$C_n\subset\omega_1$ such that $h_{t_n}^{-1}((0,\gamma))$ misses the fibers above $C_n$, that is,
$h_{t_n}(\pi^{-1}(C_n))\subset [\gamma,\omega_1)$.
By continuity, $h_t(\pi^{-1}(B))\subset [\gamma,\omega_1)$ for $B=C\cap\cap_{n\in\omega}C_n$, a contradiction.\\
\endproof

We believe that the following problem has a negative answer, but cannot back our claim with any evidence.

\begin{prob}\label{prob:bundles}
  Is it possible to have a locally trivial fiber bundle projection $p:E\to\LL_+$ with a contractible total space $E$~?
\end{prob}

%%%%%%%%%%%%%%%%%%%%%%%%%%%%%%%%%%%%%%%%%%%%%%%%%%%%%%%%
%%%%%%%%%%%%%%%%%%%%%%%%%%%%%%%%%%%%%%%%%%%%%%%%%%%%%%%%

\section{A non-contractible sub-surface of the Pr\"ufer surface}\label{sec:prufer}

We first recall the definition of the Pr\"ufer surface. We just sketch the construction and refer
to \cite[Ex. I.25--28]{GauldBook} for more details. Notice however that our notation is slightly different.
Let $H$ be the half plane $\{\langle x,y\rangle\in\R^2\,:\,x\ge 0\}$, and $H_0$ contain the points of $H$ with horizontal
coordinate $>0$.
By {\em Pr\"uferizing} $H$ at height $c\in\R$, we mean deleting $\langle 0,c\rangle\in H$
and replacing it by a copy $\R_c$ of $\R$, with the following topology.
Given $E\subset\R$, we write $E_c$ for the copy of $E$ in $\R_c$.
The neighborhoods of the non-deleted points in $H$ 
are unchanged, and 
a neighborhood basis of $u\in\R_c$ is given by all
$$T_{a,b,\epsilon}^c = (a,b)_c\cup\{ \langle x,y\rangle\in H_0 \,:\, ax+c < y < bx+c \,,\, x<\epsilon\},$$
for
$a,b\in\R_c$, $a<u<b$ and $\epsilon >0$. In words, $T_{a,b,\epsilon}^c\cap H_0$ is the triangular region 
between the lines of slopes $a,b$ passing through $\langle 0,c\rangle$ and the vertical line $x=\epsilon$.
If $A\subset\R$, by Pr\"uferizing at $A$ we mean Pr\"uferizing at each point of $A$.
If one Pr\"uferizes at $A$ and deletes $\{0\}\times(\R-A)$ from $H$, the
resulting space $P_A$ is a surface with boundary $\cup_{c\in A}\R_c$. 
A finite $A$ yields a surface homeomorphic to $H$ minus finitely many points 
in the vertical axis. 
The surface (with boundary) $P_\R$
obtained by Pr\"uferizing at {\em each} height is usually called the (separable) Pr\"ufer surface
and will be denoted by $\mathsf{P}$.
It is known for a long time that $P_A$ is contractible for any $A\subset\R$,
in particular $\mathsf{P}$ and its Moorized version (see below) are both contractible, see \cite{CalabiRosenlicht}.
We will show that deleting some points in the boundary of $\mathsf{P}$
yields a non-contractible submanifold.
Let $0_c$ denote the $0$ point in $\R_c\subset \mathsf{P}$ and set $Z = \cup_{c\in\R}\{0_c\}$.
Then $Z$ is a closed discrete subset of $\mathsf{P}$.

\begin{example}\label{ex:prufer1}
   The separable surface with boundary $\mathsf{P}-Z$ 
   has trivial homotopy groups but is non-contractible
   (and hence not heCWc).
\end{example}

We first show the following variation, for which we have a simpler proof of the non-contractibility. 

\begin{example}\label{ex:prufer2}
   Let $d>0$ and $N_d = \cup_{c\in\R}[-d,d]_c$.
   Then the separable surface with boundary $\mathsf{P}-N_d$ 
   has trivial homotopy groups but is non-contractible
   (and hence not heCWc).
\end{example}

\begin{figure}[h]
  \begin{center}
  \epsfig{figure=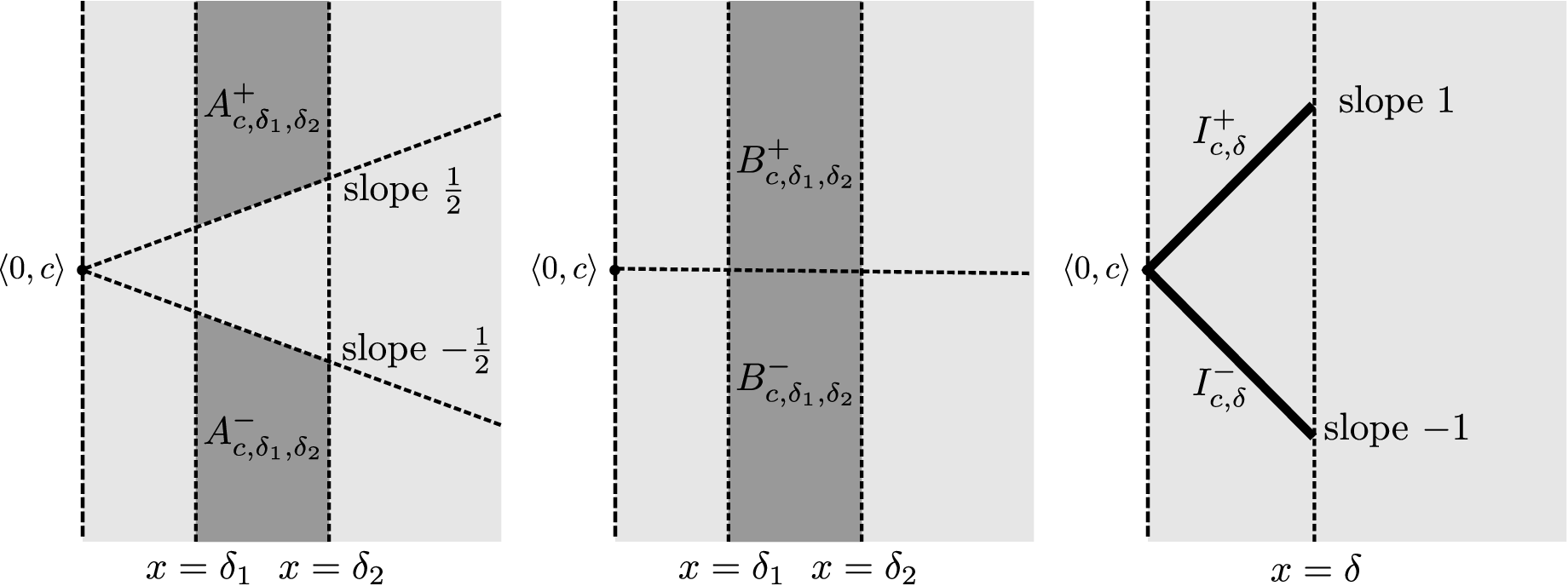, width = .8\textwidth}
  \caption{Subsets of $H_0$.}
  \label{fig:subsetsofH}
  \end{center}
\end{figure}

\proof[Details]
   Recall that $H_0$ is the open half plane 
   $\{\langle x,y\rangle\in\R^2\,:\,x > 0\}$.
   Firstly, it is clear that the homotopy groups of $\mathsf{P}-N_d$ are trivial. Indeed,
   any compact subset intersects finitely many $\R_c$, and 
   is thus contained in a submanifold homeomorphic to $H$
   minus finitely many points in its boundary, which is contractible.
   Secondly, if $d,d'\not=0$, 
   the map $\mathsf{P}-N_d \to \mathsf{P}-N_{d'}$ defined by
   $\langle x,y\rangle\mapsto \langle\frac{d'}{d}x,y\rangle$
   for points in $H_0$ and $u\mapsto \frac{d'}{d}u$ for $u\in\R_c$
   is a homeomorphism.
   We may thus assume that $d=\frac{1}{2}$ and prove that $\mathsf{P}-N_{1/2}$ is not contractible.\\
   By way of contradiction, assume $h:(\mathsf{P}-N_{1/2})\times[0,1]\to\mathsf{P}-N_{1/2}$ to be such that $h_1 = id$ and $h_0\equiv *$.
   We may assume that $*\in H_0$.  
   Write $1_c,-1_c$ for the points $1,-1\in\R_c\subset\mathsf{P}-N_{1/2}$.
   We are going to use the subsets of $H_0$ depicted in Figure \ref{fig:subsetsofH}.
   (Strictly speaking, $\langle 0,c\rangle$ belongs to $H$ and not to $H_0$, but 
   it should cause no confusion since it also does not belong to any of the subsets we consider.)
   We allow $\delta_1$ to be $0$ and $\delta_2$ to be equal to $+\infty$. 
   $I^{\pm}_{c,\delta}$ are closed segments in $H_0$, while $A^\pm_{c,\delta_1,\delta_2}, B^\pm_{c,\delta_1,\delta_2}$
   are open (hence dashed lines do not belong to them).

   Notice that in $\mathsf{P}-N_{1/2}$, for any $\delta>0$ the following hold:
   \begin{enumerate}
      \item $1_c$ is the only accumulation point of $I^+_{c,\delta}$ not in $I^+_{c,\delta}$
            and $-1_c$ is the only accumulation point of $I^-_{c,\delta}$ not in $I^-_{c,\delta}$,
      \item For any $\delta > 0$,
            $A^+_{c,0,\delta}\cup (\frac{1}{2},+\infty)_c$ is a neighborhood of $(\frac{1}{2},+\infty)_c$
            and $A^-_{c,0,\delta}\cup (-\infty,-\frac{1}{2})_c$ is a neighborhood of $(-\infty,-\frac{1}{2})_c$.
   \end{enumerate} 

\begin{claim}\label{claim:1}
   Assume that there are $\delta>0$, $u,v\in\R\cup\{\pm\infty\}$ and $t_0\in[0,1]$
   such that 
   both $h_{t_0}(1_c)$ and $h_{t_0}(-1_c)$ are in $H_0$ and have horizontal coordinate $>\delta$
   for each $c\in(u,v)$. Then:
   \\
   (i) $\exists n\ge 1,\, u^+,v^+,a^+,b^+\in\Q\cap[0,1]$, with $t_0<a^+<b^+$ and $u\le u^+<v^+\le v$
       such that both $h_t(1_c)$ and $h_t(-1_c)$ are in the closure of
       $B^+_{c,\delta/(n+1),\delta/n}$ for each $c\in[u^+,v^+]$ and $t\in[a^+,b^+]$,\\
   (ii) $\exists m\ge 1,\, u^-,v^-,a^-,b^-\in\Q\cap[0,1]$, with $t_0<a^-<b^-$ and $u\le u^-<v^-\le v$
       such that both $h_t(1_c)$ and $h_t(-1_c)$ are in the closure of
       $B^-_{c,\delta/(m+1),\delta/m}$ for each $c\in[u^-,v^-]$ and $t\in[a^-,b^-]$.
\end{claim}

\begin{figure}[h]
  \begin{center}
  \epsfig{figure=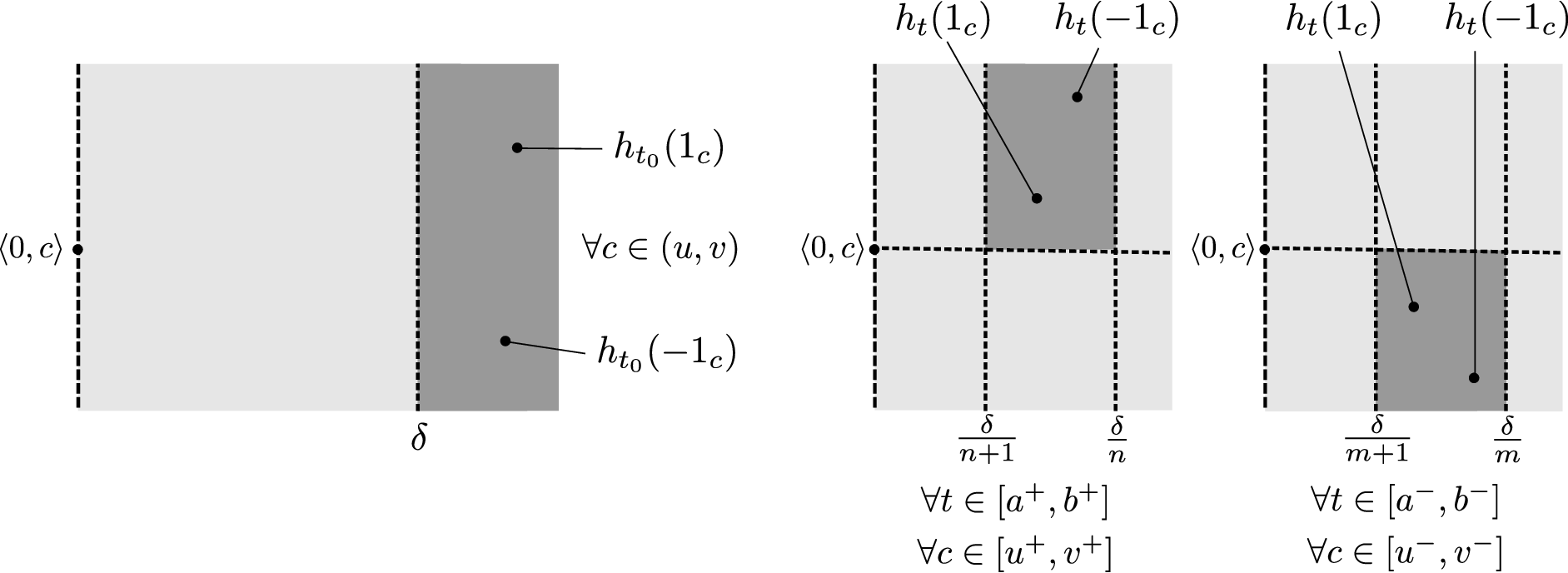, width = .95\textwidth}
  \caption{Claim \ref{claim:1}.}
  \label{fig:claim1}
  \end{center}
\end{figure}

  This claim is summarized in Figure \ref{fig:claim1}, 
  where it is to be understood that the left part implies the right part.
  Let us show why it implies the non-contractibility of $\mathsf{P}-N_{1/2}$.
  Since $h_0$ sends everything to $*\in H_0$, set $\delta_0$ to be half of $*$'s horizontal coordinate.
  We now use induction on $k$ and the claim to find nested intervals $[u_k,v_k]\subset[u_{k-1},v_{k-1}]$,
  $t_k \ge t_{k-1}$ and $\delta_k = \delta_{k-1}/(n(k)+1)$, $\delta'_k = \delta_{k-1}/(n(k))$ for some $n(k)\ge 1$ 
  such that the following holds:
  \\
  (i) $h_{t_k}(1_c)$ and $h_{t_k}(-1_c)$ are
  both in $\wb{B^+_{c,\delta_k,\delta'_k}}$ for $c\in[u_k,v_k]$ when $k$ is even,\\
  (ii) $h_{t_k}(1_c)$ and $h_{t_k}(-1_c)$ are
  both in $\wb{B^-_{c,\delta_k,\delta'_k}}$ if $c\in[u_k,v_k]$ when $k$ is odd.\\
  (These two points correspond exactly to (i) and (ii) in the claim.)
  Let $t = \sup_{k\in\omega} t_k$ and $c\in\cap_{k\in\omega}[u_k,v_k]$.
  Then $h_{t_k}(1_c)$ and $h_{t_k}(-1_c)$ should converge to $h_t(1_c)$ and $h_t(-1_c)$,
  but are divergent in $\mathsf{P}-N_{1/2}$
  since they zigzag around the horizontal at height $c$
  while going to the boundary. (This horizontal line converges to $0_c$ in $\mathsf{P}$,
  which is not present in $\mathsf{P}-N_{1/2}$.)
  \endproof

\begin{figure}[h]
  \begin{center}
     \epsfig{figure=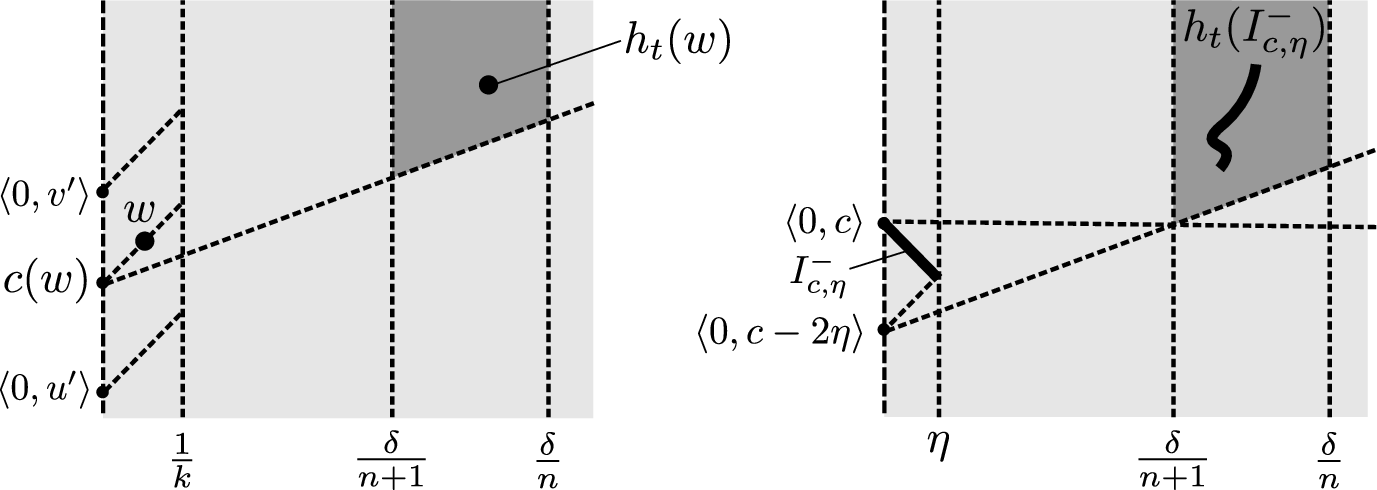, width = .7\textwidth}\\
  \caption{Proof of Claim \ref{claim:1}.}
  \label{fig:proofclaim1}
  \end{center}
\end{figure}

\proof[Proof of Claim \ref{claim:1}]
  We prove (i), as (ii) is entirely similar.
  Let $\pi:\mathsf{P}\to\R$ be the projection on the first coordinate on $H_0$
  and take constant value $0$ on the boundary components $\R_c$. Then $\pi$ is easily seen to be continuous.
  For $n,k\ge 1$ and $p,q\in\Q\cap[0,1]$, Let 
  $$ C^+_{n,k,p,q} = \{ c\in (u,v)\,:\, h_t(I^+_{c,1/k})\subset A^+_{c,\delta/(n+1),\delta/n}\,\forall t\in(p,q)\}. $$
  We show that each $c\in(u,v)$ belongs to some $C^+_{n,k,p,q}$.
  By assumption $\pi\circ h_{t_0}(1_c)$ and $\pi\circ h_{t_0}(-1_c)$ are $>\delta$.
  Notice that $(\frac{1}{2},+\infty)_c$ is a connected component of $(\mathsf{P}-N_{1/2})-H_0$,
  hence $1_c$ must first enter $H_0$ before potentially entering any other boundary component,
  and $1_c$ {\em does} leave $(\frac{1}{2},+\infty)_c$ since $\pi\circ h_{t_0}(1_c)=\delta$.
  Since $A^+_{c,0,1}\cup (\frac{1}{2},+\infty)_c$ is a neighborhood of 
  $(\frac{1}{2},+\infty)_c$, $1_c$ enters $A^+_{c,0,1}$ after leaving $(\frac{1}{2},+\infty)_c$.
  There is thus an interval $(p,q)$ (with $p>t_0$) and $n\ge 1$ such that
  $h_t(1_c)\in A^+_{c,\delta/(n+1),\delta/n}$ for each $t\in(p,q)$.
  Now, $\{1_c\}\cup I^+_{c,\eta}$ being compact for each $\eta>0$,  
  if we shrink the interval $(p,q)$ and take $\eta$ small enough,
  then $h_t(I^+_{c,\eta})\subset A^+_{c,\delta/(n+1),\delta/n}$ for each $t\in(p,q)$.
  Taking $k$ such that $\frac{1}{k}<\eta$, we obtain that $c\in C^+_{n,k,p,q}$.
  \\
  Since any interval of $\R$ is a Baire space and there are countably many $C^+_{n,k,p,q}$,
  one of them must be non-meagre. In particular, for some $n,k,p,q$, $C^+_{n,k,p,q}$ is dense
  in an interval $[u',v']\subset (u,v)$.
  If $w=\langle x,y\rangle\in H_0$, set $c(w) = y-x$.
  By continuity of $h$, if $w$ is in the parrallelogram in $H_0$
  defined by the lines of slope $1$ emerging from $\langle 0,u'\rangle$ and $\langle 0,v'\rangle$,
  the vertical line at $\frac{1}{k}$ and the boundary, then $w$ is sent to 
  the closure of $A^+_{c(w),\delta/(n+1),\delta/n}$
  by $h_t$ for $t\in(p,q)$,
  as shown on the lefthandside of Figure \ref{fig:proofclaim1}.

  By looking at righthandside of the same figure, we see that if we take $c\in (u',v')$ and $\eta$
  small enough,
  then $I^-_{c,\eta}$ is sent in the closure of $B^+_{c,\delta_{n+1},\delta/n}$ by $h_t$ for $t\in(p,q)$.
  (For those keen on deciphering formulas, 
  take $\eta$ such that $u'-2\eta > v$ and $\eta < \min(\frac{\delta}{4(n+1)},\frac{1}{k})$.) 
  Since $-1_c$ is in the closure of $I^-_{c,\eta}$, this holds also for $-1_c$. 
  To finish the proof, take a closed subinterval $[u^+,v^+]\subset (u',v')$.
\endproof

The proof that $\mathsf{P}-Z$ is also non-contractible (Example \ref{ex:prufer1})
is very similar, but with extra steps. A shortcut would be to answer the
following question in the affirmative, but we could not decide whether 
it is the case or not.

\begin{q} 
   Given $d>0$,
   are $\mathsf{P}-N_d$ and $\mathsf{P}-Z$  homeomorphic ?
\end{q}

Let us summarize the argument of $\mathsf{P}-Z$'s non-contractibility.
Assume as before that $h:(\mathsf{P}-Z)\times[0,1]\to\mathsf{P}-Z$ is a contraction with $h_1=id$ and $h_0\equiv *$.
We also consider $1_c,-1_c$.
Set $t^+(c) = \sup\{t\in[0,1]\,:\,h_t(1_c)\in (0,+\infty)_c\}$.
Then $0<t^+(c)$ and $h_{t^+(c)}(1_c)$ is a member of $(0,+\infty)_c$ that we denote by $s(c)$.
Define $A^{\pm,\ell}_{c,\delta_1,\delta_2}$
as $A^\pm_{c,\delta_1,\delta_2}$, but with slope $1/\ell$ instead of $1/2$.
Take $\ell$ such that $\frac{1}{\ell} < \frac{s(c)}{2}$.
Then there are some $p,q\in\Q$, $p<q<1$, and $n\in\omega$ such that $h_t(1_c)$ is
contained in $A^{+,\ell}_{c,0,+\infty}$ if $p<t$ and in $A^{+,\ell}_{c,1/(n+1),1/n}$
if $p<t<q$. Up to shrinking $(p,q)$,
the same holds for $I^+_{c,\delta}$ if $\delta$ is small enough.
The same Baire argument as before
yields an interval $[u,v]$ and $\ell,n,p,q$
such that for each $c\in[u,v]$,
$h_t(1_c)\in A^{+,\ell}_{c,1/(n+1),1/n}$ if $p<t<q$.
Then $h_t(I^-_{c,\eta})\subset A^{+,\ell}_{c,1/(n+1),1/n}$ as well for the same $t$'s
and $\eta$ small enough (depending on $\ell$). This shows in particular that $h_t(-1_c)$
is in the closure of $B^+_{c,1/(n+1),1/n}$ for the same $t$'s.
Now proceed with the $-1_c$, etc, to find the same zigzaging diverging
sequence.

\vskip .3cm
Our examples so far all have a boundary, but we may get rid of it by Mooreizing the boundary components
(see again \cite[Ex. I.27]{GauldBook} for details).
In its classical version,
we take a boundary component $\R_c$ obtained by Pr\"uferization and identify $x$ with $-x$ in it.
This has the effect of folding $\R_c$ which is thus `swallowed' by the manifold and is not a boundary component anymore.
We can perform this identification in all the boundary components or only in some of them.
If one takes the Pr\"ufer surface and Moorize each boundary component, we obtain a contractible boundaryless 
non-metrizable surface.
We can also fold around another point than $0_c$,
or use more complicated functions for the identification. 
For instance, if as above $0_c$ has been removed from the boundary component $\R_c$,
we may fold $(0,+\infty)_c$ around $1_c$ by using any homeomorphism between $(0,1]_c$ and $[1,+\infty)_c$ for the
identification, and similarly fold $(-\infty,0)_c$ around $-1_c$.

\begin{example}
    The separable surface obtained by Mooreizing the boundary components of $\mathsf{P}-Z$,
    by folding $(0,+\infty)_c$ around $1_c$ and $(-\infty,0)_c$ around $-1_c$
    is non-contractible.   
\end{example}
The proof of non-contractibility in Example \ref{ex:prufer1} works verbatim.
Another way to get rid of boundary components is to glue collars to them, that is,
copies of $\R\times[0,1)$. Doing it to $\mathsf{P}$ yields a contractible non-separable surface, while doing it to 
$\mathsf{P}-Z$ yields a non-contractible subsurface of the former. For more on these variations of the Pr\"ufer surface,
see \cite{MesziguesBridges}.

%%%%%%%%%%%%%%%%%%%%%%%%%%%%%%%%%%%%%%%%%%%%%%%%%%%%%%%%%%%%%%%%%%

\section{$\omega_1$-trees and their road spaces}\label{sec:trees}

The contents of this section appeared in the preprint \cite{meszigues-contract-trees}.
The objects we study are set theoretic trees and their road spaces, 
which are obtained by joining consecutive points by a line segment with a topology that makes it locally 
embeddable in $\R^2$ (the details are given below).
These spaces are good toy models for Type I non-metrizable manifolds.

Recall that
a {\em tree} $T$ is a partially ordered set such that each point has a well ordered set of predecessors.
We define the {\em height} of $x\in T$ and of $T$, the $\alpha$-th-level $\lev_\alpha(T)$, the {\em chains} and {\em antichains} in $T$ as 
usual, see for instance \cite[Section II.5]{kunen} or \cite{Nyikos:trees}. A subset $E$ of $T$ is {\em order-dense}
iff for any $y\in T$ there is $y\in E$ with $x<y$.
All trees in this note are endowed with the order (also called interval) topology. 
A tree is {\em rooted} iff it has a unique minimal element called the {\em root}
and {\em binary} iff each member has exactly $2$ immediate successors.
A {\em subtree} $S$ is a subset of $T$ with order restricted to $S$. 
Notice that the induced topology on a subtree $S$ is finer (sometimes strictly)
than the one given by the 
order restricted to $S$. Both topologies agree if $S$ is closed in $T$.
We assume that our trees are Hausdorff, 
that is, if $x,y\in T$ are at a limit level and have the same predecessors, then $x=y$. (This could be false for a subtree.)
An {\em $\omega_1$-tree} has countable levels and height $\omega_1$.
A tree is {\em Aronszajn} (resp. {\em Suslin}) if it has height $\omega_1$ and its chains 
(resp. chains and antichains) are at most countable.
When $x\in T$ and $\alpha$ is an ordinal, write $T_{\ge x} = \{y\in T\,:\, y\ge x\}$,
$x\upharpoonright\alpha$ for the unique predecessor of $x$ at level $\alpha$ (if $x$ is below the $\alpha$-th level, $x\upharpoonright\alpha=x$)
and $T_{\le\alpha}$ for the subset of elements at level $\le\alpha$. 
If $E\subset T$, set $E^\downarrow =\{x\in T\,:\,\exists y\in E\text{ with }x\le y\}$ to be its downward closure.
We say that the tree $T$ is {\em $\R$-special} iff there is a strictly increasing (not necessarily continuous)
function $T\to \R$. Recall that $\R$-special binary $\omega_1$-trees exist in {\bf ZFC}.

\begin{figure}[h]
  \begin{center}
  \epsfig{figure=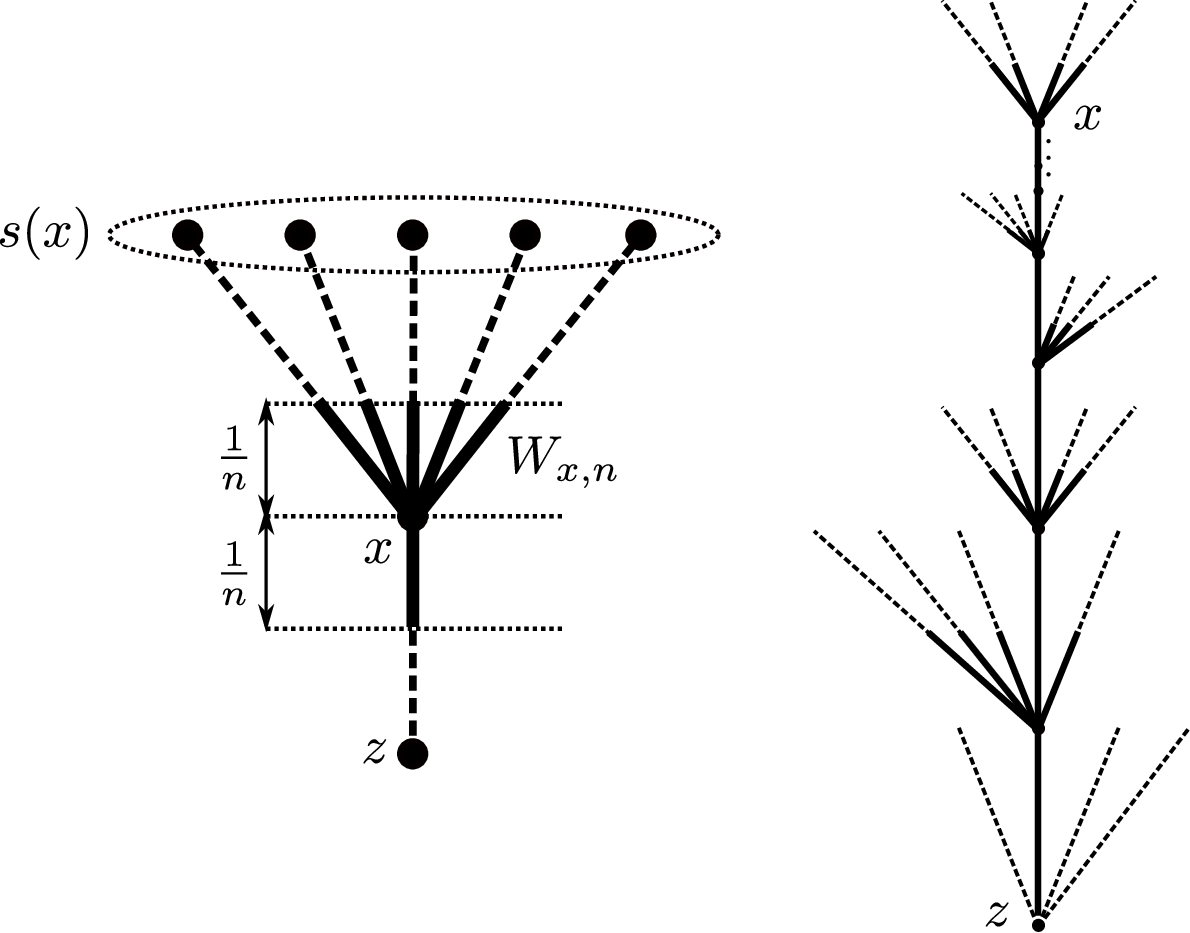,width=.5\textwidth}
  \caption{Neighboroods in $R_T$.}
  \label{fig:roadspace}
  \end{center}
\end{figure}

The road space $R_T$ of a tree $T$ is obtained by joining consecutive points by an interval $[0,1]$, 
with $0$ glued to the lowest point and $1$ to the highest. We extend the order in the obvious way.
When convenient we consider $T$ as a subset of $R_T$.
The topology on the interior of the added intervals is that of $(0,1)$.
For $x\in T\subset R_T$, in order for 
$R_T$ to be (locally) connected any open set containing $x$ must contain a small portion of each interval emanating from $x$. 
In order for the space to be locally metrizable (and hence first countable),
we take these portions uniformly as follows.
Denote by $[0,1]_x^y\subset R_T$
the interval between the two consecutive points $x,y\in T$. If $A\subset[0,1]$, then $A_x^y$ is understood as the 
corresponding subset of $[0,1]_x^y$.
For singletons we usually write $a_x^y$ instead of $\{a\}_x^y$.
If $x\in T$, denote by $s(x)$ the set containing its immediate successors and set 
$W_{x,n} =\displaystyle \bigcup_{y\in s(x)} [0,1/n)_x^y$.
If $x\in s(z)$ with $z\in T$, a basis for the neighborhoods of $x$ in $R_T$ is given by 
$\{ W_{x,n} \cup (1-1/n,1]_z^x\,:\,n\in\omega\}$  (see Figure \ref{fig:roadspace}, left).
A basic neighborhood of $x\in T\subset\R_T$ at a limit level
is obtained by 
choosing some $z\in T$, $z<x$ and $n\in\omega$ and taking the segment $\{y\in R_T\,:\,z<y<x\}$ union each $W_{w,n}$ 
for those $w\in T$ with $z<w\le x$ (see Figure \ref{fig:roadspace}, right).
This makes $R_T$ locally embeddable in $\R^2$ (as seen by induction, or by extrapolating from
Figure \ref{fig:roadspace}).
The induced topology on $T\subset R_T$ is that of $T$, and $R_T$ is arc connected if and only if $T$ is rooted.

\vskip .3cm
The first result of this section is the following.

\begin{thm}\label{thm:main1}
   The road space $R_T$ of a rooted $\R$-special tree is contractible.
\end{thm}

This gives half of the proof of Theorem \ref{thm:main-contractible}, and
is proved in Subsection \ref{subsec:contractible}.
The other two theorems of this section are generalizations of Theorem \ref{thmdeogauld}.
We depart somewhat from the manifold theme of this paper and investigate quite precise
properties of subspaces, trying to see whether they prevent global contractibility.
The first version is to weaken the assumption in Theorem \ref{thmdeogauld}
that the subspace is a copy of $\omega_1$.

\begin{thm}\label{thm:main2}
   Let $T$ be a tree, $S\subset T$ be uncountable and $\omega_1$-compact in the subspace topology, 
   and let $X$ be a locally metrizable space.
   Then there is no continuous $h\colon S\times[0,1]\to X$ 
   such that $h_1\colon S\to X$ has uncountable image and $h_0$ is constant.
   In particular, if $X$ contains a copy of $S$ then $X$ is not contractible.
\end{thm}

Our 
proof of Theorem \ref{thm:main2} is entirely topological, for one using forcing see \cite{meszigues-contract-trees}.
The first step is to show that one can assume without loss of generality that $T$ is an Aronszajn tree.
Then,
the bulk of our argument consists of showing
that many well known properties of Suslin trees also hold for $\omega_1$-compact subsets of trees of height $\omega_1$.
These generalizations are probably part of the folklore and no more than standard exercices, but we found
convenient to gather the proofs of most of them in Subsection \ref{subsec:trees}.
Then the `real' proof of this theorem 
(and the fact that it implies Theorem \ref{thm:main-w1cpct}) is contained in Subsection \ref{subsec:w1cpct}.
The property that is central to our argument is that any continuous map from $S$ to a locally metrizable space
must be constant on the points above $x\in S$ if $x$ is high enough in the tree (see Lemma \ref{lemgen:1} (i) below).
This is a slight generalization of a result of Stepr\=ans \cite{StepransTrees}.\\
Notice that an uncountable $\omega_1$-compact subset of a tree cannot be metrizable, since $\omega_1$-compactness
is equivalent to Lindel\"ofness in metrizable spaces (see, e.g., \cite[Thm 4.1.15]{Engelking}), 
and a Lindel\"of subset of a tree with countable levels is countable.
\\
A quick corollary, which is another exhibit of 
the similarities between $\omega_1$ and Suslin trees, is the following.
Its proof is immediate since the road space of a tree satisfying the assumptions is locally metrizable and
a Suslin tree is $\omega_1$-compact, as well known (see Lemma \ref{lemgen:0} below).

\begin{cor}\label{cor:main1}
   Let $T$ be a tree of height $\omega_1$. 
   If $T$ contains an uncountable subset $S$ which is $\omega_1$-compact in the subspace topology
   (in particular, if $T$ is a Suslin tree),
   then its road space $R_T$ is non-contractible.
\end{cor}

Our second generalization of Theorem \ref{thmdeogauld} is to weaken the local metrizability of the target space $X$. 
First countability alone is not enough, as shown by Example \ref{ex:cones}.
What we were able to prove is the following.

\begin{thm}\label{thm:main3}
  Let $X$ be a regular space with $G_\delta$ points such that each point has a closed $[\aleph_0,\aleph_1]$-compact
  neighborhood.
  Let $S$ be a stationary subset of $\omega_1$ endowed with the subspace topology. 
  Then there is no continuous $h\colon S\times[0,1]\to X$ such that $h_1\colon S\to X$ 
  has an uncountable image and $h_0$ is constant.
  In particular, if $X$ contains a copy of $S$ then $X$ is not contractible.
\end{thm}
Recall that a space is {\em $[\aleph_0,\aleph_1]$-compact} if and only if any
cover of $X$ of size $\le\aleph_1$ has a countable subcover. This is equivalent to
the property that given any uncountable subspace $E$ of $X$ there is $x\in X$ which is a {\em condensation point of $E$}, that is,
any neighborhood of $x$ contains uncountably many points of $E$.   
We do not know whether this theorem holds if $S$ is a Suslin tree but note in passing that 
our main tool (any continuous map from $S$ to $X$
must be constant on the points above $x\in S$ if $x$ is high enough, see Lemma \ref{lemma:omega_1} below) does not.
(This should appear in a later paper.)
Theorem \ref{thm:main3} holds if $X$ is locally metrizable.
Indeed, a stationary subset of $\omega_1$ is $\omega_1$-compact,
$\omega_1$-compactness is preserved by continuous functions and
is equivalent to Lindel\"ofness in metrizable spaces; hence the image of $S\times[0,1]$ in $X$ is (hereditarily) Lindel\"of.

Corollary \ref{cor:main1} inspires the following question. 

\begin{q}\label{q:1}
  Let $T$ be a tree and $R_T$ be its road space.\\
  (a) Does the contractibility of $R_T$ implies that $T$ is $\R$-special~?\\
  (b) Does the non-contractibility of $R_T$ imply that $T$ contains a Suslin subtree~?
\end{q}

Notice that it is not possible for (a) and (b) to both have a positive answer since
there are models of set theory with a non-$\R$-special tree $T$ 
that does not contain a Suslin subtree, see \cite{Schlindwein:2003}. We do not know whether these trees have contractible road spaces.

%%%%%%%%%%%%%%%%%%%%%%%

\subsection{Proof of Theorem \ref{thm:main1}}\label{subsec:contractible}

\proof
Let $T$ be a rooted $\R$-special tree 
with root $r$ and let $f:T\to\R$ be a strictly increasing function. By replacing $f(x)$ by $\sup_{y<x}f(y)$ 
when $x$ is at limit levels, we may assume that $f$ is continuous. By composing with a strictly increasing function, we
may assume that the range of $f$ is contained in $[0,1]$.
We first start to define the contraction as a map $h:T\times[0,1]\to R_T$.
First, set $h_1(x) = x$ and 
then define $h$ such that $x$ travels downwards in $R_T$, starting to move exactly at time $t=f(x)$ and reaching
$y<x$ ($y\in T$) exactly at time $t=f(y)$. Since $f$ is strictly increasing, there is time available to cross the interval between
consecutive points. For those who like
(less readable) formulas, let $x$ at level $\alpha$ be given.
If $f(x)\le t$, we set $h_t(x) = x$.
If $\beta<\alpha$ and $f(x\upharpoonright \beta)\le t\le f(x\upharpoonright \beta+1)$,
we set $h_t(x)$ to be $u_{x\upharpoonright \beta}^{x\upharpoonright \beta+1}$ where
   $$ u=
    \frac{t-f(x\upharpoonright \beta )}{f(x\upharpoonright \beta+1) - f(x\upharpoonright \beta)}.$$
Finally, if $t\le f(r)$, we set $h_t(x) = r$.
It should be clear that $h$ is continuous (on $T\times[0,1]$) and that $h_0(x) = r$ for all $x\in T$.
Notice that $h$ has the property that at time $f(x)$, all of $T_{\ge x}$ is squished onto $x$.
Actually, each point on the tree starts to move exactly when all the points above it reach it, all at the same time.
This enables to extend easily the map to all of $R_T$, a point in the segment $[0,1]_{x}^y$ does not move until
$y$ reaches it, and then it follows it until the end.
This gives the required contraction.
\endproof

Notice that if one sets $j_t(x) = h_{1-t}(x)$ for $t\in[0,1]$, $j_t(x) = x$ when $t<0$
and $j_t(x) = r$ when $t>1$, then $j$ is actually a flow, that is $j_t\circ j_s = j_{t+s}$.

%%%%%%%%%%%%%%%%%%%%%%%%%%%%%%%%%%%%%%%%%%%%%%%

\subsection{A collection of facts on uncountable $\omega_1$-compact subsets of trees}\label{subsec:trees}

Notations: Given a tree $T$, if $A,B\subset T$ and $x\in T$, 
$A<B$ means that 
each member of $A$ is $<$ each member of $B$. We denote $\{x\} < A$ by $x<A$
and $T_{\ge x}\cap A$ by $A_{\ge x}$.
Notice that $\text{Lev}_\alpha(A^\downarrow) = \text{Lev}_\alpha(T) \cap A^\downarrow$.

\begin{lemma}\label{lemgen:0}
   Let $T$ be a tree of height $\omega_1$ and $S\subset T$ be uncountable and endowed with the subspace topology. Then the following hold.\\
   (a) The subspace topology on $S^\downarrow$ always agrees with the topology given by the induced order.
       If $S$ is closed in $T$, then the subspace topology agrees with the topology given by the induced order on $S$.
       \\
   (b) If $S$ is $\omega_1$-compact, then it intersects a stationary subset of levels of $S^\downarrow$ (and of $T$).\\   
   (c) An antichain is closed discrete in $T$, and a closed discrete subset of $T$ is an at most countable union of antichains.\\
   (d) If $S$ is $\omega_1$-compact then so is $S^\downarrow$.\\
   (e) If $S$ is $\omega_1$-compact and does not contain an uncountable chain, then $S^\downarrow$ is Suslin.\\
   (f) If $A\subset S\subset T$, then $A$ is a maximal antichain in $S$ if and only if it is a maximal antichain in $S^\downarrow$.
\end{lemma}
\proof
(a) and (f) are straightforward. For (b), if $S$ avoids a club set of levels of $S^\downarrow$ 
then taking one point between consecutive avoided levels (when available) yields an uncountable discrete subset which is closed in $S$.
Item (c) is proved e.g. in \cite[Thm 4.11]{Nyikos:trees}, and (d)--(e) follow immediately from it.
\endproof

\begin{lemma}\label{lemgen:1}
   Let $T$ be a tree of height $\omega_1$ and $S\subset T$ be uncountable and $\omega_1$-compact in the subspace topology.
   Then the following hold.\\
   (a) $S^\downarrow$ is the disjoint union of an at most countable set, a Suslin tree and at most countably many copies of $\omega_1$.
       In particular, $S^\downarrow$ is an $\omega_1$-tree. \\
   (b) There is $\alpha\in\omega_1$ such that $|S_{\ge x}|=|(S^\downarrow)_{\ge x}| = \aleph_1$ when $x$ is above level $\alpha$.\\
   (c) If $C\subset S^\downarrow$ is uncountable, there is $x\in S$ such that 
       $(S^\downarrow)_{\ge x}\subset C^\downarrow$ and $|S_{\ge x}|=|(S^\downarrow)_{\ge x}| = \aleph_1$.\\
   (d) $C^\downarrow$ is $\omega_1$-compact for any $C\subset S^\downarrow$.\\
   (e) If $C\subset S^\downarrow$ is club in $S^\downarrow$, then 
       $\{\gamma\,:\,\text{Lev}_\gamma(C^\downarrow)\subset C\}$ is club in $\omega_1$.\\
   (f) If $C\subset S^\downarrow$ is club in $S^\downarrow$, 
       then $S\cap C$ intersects a stationary subset of levels.\\
   (g) If $E_n\subset S^\downarrow$ are closed and order-dense in $S^\downarrow$ 
       for $n\in\omega$, then $\cap_{n\in\omega}E_n$ is also closed and order-dense.
       \\
   (h) Let $E,F\subset S\subset T$ be closed in $S$.
       If $|E\cap F| \le\aleph_0$, then $|E^\downarrow\cap F^\downarrow|\le\aleph_0$.\\
   (i) Let $f:S\to Y$ be continuous where $Y$ is a metrizable space. Then there is $\beta\in\omega_1$
       such that $f(S_{\ge x})$ is a singleton whenever $x$ is above level $\beta$ in $T$.
       In particular, $f$ has a countable image.
\end{lemma}
\proof
We will use several times (without acknowledging it explicitly) the fact that an at most countable intersection of club
subsets of $\omega_1$ is club.

(a) If $S$ contains an uncountable branch $B$, then there is some minimal $x_B\in B$ such that $S_{\ge x_B}$ is linearly ordered.
       Indeed, otherwise $S$ contains an uncountable antichain (take points branching away from $B$ above each height)
       and hence an uncountable closed discrete subset by Lemma \ref{lemgen:0} (c).
       Since the minimal elements of $\{x\,:\,S_{\ge x}\subset S\text{ is an uncountable branch}\}$ is an antichain, 
       $S$ contains at most countably many disjoint uncountable branches.
       Notice that a maximal branch in $S^\downarrow$ contains an unbounded branch of $S$.
       Removing the branches above each minimal
       $x_B$ in $S^\downarrow$, what remains is either countable or a Suslin tree by Lemma \ref{lemgen:0} (e).

(b) follows immediately from (a) and the equivalent statement for Suslin trees.

(c) If $C\cap B$ is uncountable for some branch $B\subset S^\downarrow$, we are over. 
     If not, by (a) and (b)
     we can assume that $S^\downarrow$ is a Suslin tree such that
     $|S_{\ge x}|=\aleph_1$ for each $x\in S^\downarrow$. 
     Then the result is well known (see, for instance, the claim in Theorem 2.1 in \cite{DevlinShelah}).

(d) is immediate by (a) since an uncountable (downward closed) subset of a Suslin tree is a Suslin tree.

(e) By (d), $C^\downarrow$ is $\omega_1$-compact. Fix $\alpha$ given by (b) and some $\beta>\alpha$.
    Since each uncountable maximal branch $B$ of $C^\downarrow$ is a copy of $\omega_1$, $C$ contains a club set 
    of levels of $B$. By (a) we may assume that $C^\downarrow$ is Suslin.
    By Lemma \ref{lemgen:0} (f), 
    for each $n\in\omega$ we may find a countable antichain $A_n\subset C$ above level $\beta$
    which is maximal in $C^\downarrow$ and such that 
    $A_{n+1}>A_n$.
    Let $\gamma = \sup_{n\in\omega}\sup\{\text{height}(x)\,:\,x\in A_n\}$.
    By construction $\text{Lev}_\gamma(C^\downarrow)$ is the set of limit points of $\cup_{n\in\omega}A_n$, hence
    since $C$ is closed $\text{Lev}_\gamma(C^\downarrow)\subset C$.
    This shows that $\{\gamma\,:\,\text{Lev}_\gamma(C^\downarrow)\subset C\}$ is unbounded in $\omega_1$,
    and closedness is obvious.

(f) If $C\cap B$ is unbounded for some maximal uncountable branch $B\subset S^\downarrow$, then 
       it is homeomorphic to a club subset of $\omega_1$.
       By Lemma \ref{lemgen:0} (b) $S\cap C \cap B$ is thus stationary. 
       If $C\cap B$ is bounded for each uncountable branch of $S^\downarrow$,
       we may assume by (a) that $S^\downarrow$ is a Suslin tree in which $C$ is unbounded.
       By (c) it follows that $S$ is unbounded in $C^\downarrow$ as well. 
       By (e) and Lemma \ref{lemgen:0} (b),
       $S$ intersects $C$ on a stationary set of levels.

(g) Closedness is immediate, and order-density follows immediately by (a) and the fact that 
    the result holds for Suslin trees and $\omega_1$.
\iffalse
       hence let us show that $\cap_{n\in\omega}E_n$ is order-dense.
       First, since each uncountable maximal branch $B$ of $S^\downarrow$ is homeomorphic to $\omega_1$,
       $\cap_{n\in\omega_1} E_n$ is club in $B$. By (a) we may now assume that $S^\downarrow$ is Suslin.
       Let $\alpha$ be given by Lemma \ref{lemgen:1} (a).
       Let $x\in S^\downarrow$. We may assume that $x$ is at level above $\alpha$.
       Fix $\sigma:\omega\to\omega$ such that $|\sigma^{-1}(\{n\})| = \aleph_0$ for each $n\in\omega$.
       Since each $E_n$ is order-dense, for each $n$ we may take an antichain $A_n\subset E_{\sigma(n)}$ which is maximal in $S^\downarrow$
       such that $x<A_n < A_{n+1}$. Since $S^\downarrow$ is Suslin there is some $y\in T$ which is the limit of a sequence $x_n>x$
       such that $x_n\in A_n$. It follows that $y\in\cap_{n\in\omega}E_n$.
\fi

(h) 
   Let $\wb{E},\wb{F}$ be the closures in $S^\downarrow$ of $E,F$.
   If $\wb{E}\cap\wb{F}$ is unbounded in $S^\downarrow$, by (f) $S\cap\wb{E}\cap\wb{F} = (S\cap\wb{E})\cap (S\cap\wb{F}) = E\cap F$
   is unbounded, a contradiction.
   Hence $\wb{E}\cap\wb{F}$ is bounded and thus $\wb{E}$ and $\wb{F}$ are disjoint above some level $\alpha$.
   It follows that $E$ and $F$ cannot be both unbounded in the same uncountable branch. 
   It is well known (see e.g. \cite[Thm 6.18]{Todorcevic:1984} or \cite[Thm 2.1]{DevlinShelah}) 
   that if $A,B$ are disjoint closed sets in a Suslin tree,
   then $A^\downarrow\cap B^\downarrow$ is at most countable. 
   Together with (a), this shows that $|E^\downarrow\cap F^\downarrow|\le\aleph_0$.
   
   (i) Our proof is a slight adaptation of Stepr\=ans topological proof in \cite{StepransTrees} 
       that a real valued map with domain a Suslin tree has a countable image. 
       Denote the distance in $Y$ by $\text{dist}(\cdot,\cdot)$.
       By (b) me may assume that $S_{\ge x}$ is uncountable for each $x$.
       Set $D(\epsilon) = \{ x\in S\,:\,\text{diam}(f (S_{\ge x})) \le \epsilon\}$, where diam stands for diameter, that is, the supremum
       of the distances between two points in a set.
       Assume for now that $D(\epsilon)$ is order-dense in $S$ (and thus in $S^\downarrow$) when $\epsilon >0$.
       Let $\wb{D(\epsilon)}$ denote the closure of $D(\epsilon)$ in $S^\downarrow$.
       By (g), $D = \cap_{n\in\omega,n>0}\wb{D(1/n)}$ is closed and order dense in $S^\downarrow$, hence
       by (f) $S\cap D$ intersects stationary many levels above each $x\in S$ (since $|S_{\ge x}|=\aleph_1$).
       Moreover, $D$ is upward closed in $S$. 
       Denote by $A$ the minimal elements of $D\cap S$. Then $A$ is an antichain of $S^\downarrow$,
       let $\beta$ be the supremum of the heights of its members.
       For each $x\in S$ above level $\beta$, the diameter of $f\bigl(S_{\ge x}\bigr)$ is $0$, hence  
       $f\bigl(S_{\ge x}\bigr)= f(\{x\})$ and the lemma is proved.\\
       To finish, we now prove that $D(\epsilon)$ is order-dense in $S$.
       Suppose that it is not the case and let $x\in S$ be such that 
       $\text{diam}\bigl(f(S_{\ge y})\bigr) >\epsilon$ for each $y>x$, $y\in S$.
       We build inductively antichains $A^\alpha_n$ ($n\in\omega,\alpha\in\omega_1$) such that the following hold.\\
         $\bullet$ $A_n^\alpha\subset S$ is maximal above $x$, that is, in $(S^\downarrow)_{\ge x}$,\\
         $\bullet$ $A_{n+1}^\alpha > A_n^\alpha > A_m^\beta$ for each $n,m\in\omega$ and $\alpha > \beta$,\\
         $\bullet$  If $u\in A_n^\alpha$, $v\in A_{n+1}^\alpha$, then $\text{dist}(f(u),f(v))\ge\epsilon/4$.\\
       Assume that $A_n^\alpha$ is defined.
       Set 
         $$E=\{ z\in S\,:\, z>A_n^\alpha\text{ and dist}(f(z),f(u))\ge\epsilon/4 
           \text{, where $u$ is the member of $A_n^\alpha$ below $z$}\}.$$ 
       It is enough to see that $E$ is order-dense in $(S^\downarrow)_{\ge x}$, 
       since then we may put its minimal elements in $A_{n+1}^\alpha$.
       Let thus $w\in (S^\downarrow)_{\ge x}$, $w>u\in A_n^\alpha$.
       Up to going further up, we may assume that $w\in S$.
       If $\text{dist}(f(w),f(u))\ge\epsilon/4$, then $w\in E$.
       If not, then 
       $\text{dist}(f(w),f(u))<\epsilon/4$.
       Choose $v\in S_{\ge w}$ such that $\text{dist}(f(w),f(v))>\epsilon/2$ 
       (which exists since we assumed $\text{diam}\bigl(f(S_{\ge y})\bigr) >\epsilon$ for each $y>x$).
       Then 
         $$\text{dist}(f(v),f(u)) \ge \text{dist}(f(v),f(w)) - \text{dist}(f(w),f(u)) >\epsilon /4,$$ 
       and $v\in E$.
       If $A_n^\gamma$ is chosen for each $n\in\omega$ and each $\gamma<\alpha$, set $A_0^\alpha$ to be an antichain
       in $S$, maximal in $(S^\downarrow)_{\ge x}$, 
       whose members are all $>\cup_{n\in\omega,\gamma<\alpha}A_n^\gamma$. This defines $A_n^\alpha$
       for each $n\in\omega,\alpha\in\omega_1$.\\
       Set $\beta(\alpha)$ to be $\sup\{ \text{height}(y)\,:\,y\in\cup_{n\in\omega} A_n^\alpha\}$,
       let $C$ be the closure in $\omega_1$ of $\{\beta(\alpha)\,:\,\alpha\in\omega_1\}$,
       and $C'$ be its derived set (that is its limit points).
       By construction, if there is $y\in S_{\ge x}$ whose height in $S^\downarrow$ is in $C'$, then $f$ is not continuous at $y$
       as there is a sequence of points in $S$ converging to $y$ whose images are $\ge\epsilon/4$ apart.
       But by Lemma \ref{lemgen:0} (b) (and the fact that $S_{\ge x}$ is $\omega_1$-compact), there must be such an $y$, a contradiction. 
       This shows that $D(\epsilon)$ is order-dense and concludes the proof.
\endproof

%%%%%%%%%%%%%%%%%%%%%%%%%%%%%%%%

\subsection{Proof of Theorem \ref{thm:main2}}\label{subsec:w1cpct}. 

Our proof relies on simple consequences of the properties given in Lemma \ref{lemgen:1}, especially (i).
When $S$ is a subset of a tree $T$, the height of a point of $S$ is to be understood as its height in $T$.
Recall that $S_{\ge x} = T_{\ge x}\cap S$.

\begin{lemma}
   Let $X$ be a space containing a closed metrizable $G_\delta$ subset $B\subset X$, 
   $S$ be an uncountable $\omega_1$-compact subset of a tree of height $\omega_1$ and
   $f:S\to X$ be continuous.
   Then there is $\alpha\in\omega_1$ such that
   either $f(S_{\ge x})\cap B = \varnothing$ or $f(S_{\ge x})$ is a singleton whenever $x$ is at height $\ge\alpha$.
\end{lemma}
\proof
   By Lemma \ref{lemgen:1} (b), above some level each $S_{\ge x}$ is uncountable, we assume for simplicity that this holds for each $x\in S$.
   Let $U_n$ be open sets such that $\cap_{n\in\omega} U_n = B$.
   By Lemma \ref{lemgen:1} (h), there is $\beta\in\omega_1$ such that for each $n$ we have
   $$ \left(\,\left(f^{-1}(B))\right)^\downarrow\cap \left(f^{-1}(X-U_n)\right)^\downarrow \right) \, - T_{\le\beta} = \varnothing.$$ 
   
   It follows that if $x\in S$ is at level above $\beta$,
   either $f(S_{\ge x})\subset B$ or $f(S_{\ge x})\cap B = \varnothing$.
   Let $E = \{x\in S\,:\, \lev(x)\ge\beta\text{ and }f(S_{\ge x})\subset B\}$.
   Then $E$ is an upward closed subspace of $S$, in particular it is an $\omega_1$-compact subspace of $T$.
   By Lemma \ref{lemgen:1} (i), there is $\alpha\ge\beta$ such that 
   $f$ is constant on $S_{\ge x}$ whenever $x\in E$ is at level $\ge\alpha$.
   This proves the lemma.
\endproof

\begin{cor}
   \label{cor:steprlocmet}
   Let $X$ be a space and $U\subset X$ be open such that
   $\wb{U}$ is contained in a metrizable open $V\subset X$.
   Let $S$ be an uncountable $\omega_1$-compact subset of a tree of height $\omega_1$.
   Let $h:S\times[0,1]\to X$ be continuous.
   Then there is $\alpha\in\omega_1$ such that for each $t\in[0,1]$ and each $x\in S$ above level $\alpha$,
   either $h_t^{-1}(U) \cap S_{\ge x} = \varnothing$, or $h_t$ is constant on $S_{\ge x}$.
\end{cor}
\proof
   Let $\{t_n\,:\,n\in\omega\}$ be a countable dense subset of $[0,1]$. Set $B=\wb{U}$.
   Since $B$ is closed in the metrizable subset $V$, it is a $G_\delta$.
   The previous lemma shows that there is some $\alpha$
   such that when $x$ is above level $\alpha$, either $h_{t_n}^{-1}(B)\cap S_{\ge x} = \varnothing$ or $h_{t_n}$ is constant on $S_{\ge x}$
   for each $n\in\omega$.
   The result follows by continuity.
\endproof

\proof[Proof of Theorem \ref{thm:main2}]
We first show that we may assume that $T$ is an Aronszajn tree.
Indeed, consider $S^\downarrow\subset T$;
if some level of $S^\downarrow$ is uncountable, then 
$S$ is not $\omega_1$-compact (by Lemma \ref{lemgen:0} (c) \& (d)).
If $S^\downarrow$ contains an uncountable branch $E$, 
then $E$ is $\omega_1$-compact in the subspace topology and thus homeomorphic to a 
stationary subset of $\omega_1$, and we may apply Theorem \ref{thm:main3} (whose proof given later
is independent from this one).\\
Let $h:S\times[0,1]\to X$ be continuous such that $h_0(x) = u_0\in X$ and $h_1:S\to X$ has uncountable image.
The set of $x\in S$ such that $S_{\ge x}$ has uncountable image under $h_1$ is uncountable and downward closed, 
hence by Lemma \ref{lemgen:1} (c) there is $x$
such that the image of $S_{\ge z}$ under $h_1$ is uncountable for each $z\ge x$.
Up to replacing $S^\downarrow$ by $(S^\downarrow)_{\ge x}$ we assume that this holds for all $z\in S^\downarrow$.
For $x\in S$, set 
$$
\tau(x)=\sup\{t\,:\,h_t\text{ is constant on }S_{\ge x}\}.
$$
Then $\tau$ is an increasing map $S\to\R$, however $\tau$ is a priori neither continuous nor strictly increasing.
We will show that there is a closed unbounded $C\subset\omega_1$ such that $\tau$ is strictly increasing on the subspace of members of $S$ at 
levels belonging to $C$. 
Such a subspace is uncountable (and thus unbounded) in $S$ by Lemma \ref{lemgen:1} (f).
Hence as a partially ordered space $S$ contains an $\R$-special tree and thus
(at least) an uncountable antichain (see for instance \cite[Thm 4.29]{Nyikos:trees}), a contradiction
with the fact that $S$ is $\omega_1$-compact.
\\
It is enough to show that for each $x\in S$, there is $\alpha$ such that $\tau(y)>\tau(x)$ whenever $y>x$ is at level $\ge\alpha$. Indeed,
since the levels of $S^\downarrow$ are countable by Lemma \ref{lemgen:1} (a) and $\tau$ is increasing, a simple induction provides $C$.
Let $x\in S$ be fixed.
By continuity, $h_{\tau(x)}$ is constant on $S_{\ge x}$ with value $u = h_{\tau(x)}(x)$ and thus $\tau(x)<1$ since $h_1(S_{\ge x})$ is uncountable. 
(While it is not needed, notice that $h$ restricted to the subspace $S_{\ge x} \times [\tau(x),1]$ 
'contracts' all of $S_{\ge x}$ to the point $u$.)
Since $X$ is locally metrizable, we may choose an open $U\ni u$ such that $B=\wb{U}$ is contained
in an open metrizable set.
Let $\alpha$ be given by 
Corollary \ref{cor:steprlocmet}. We may assume that $\alpha>\text{height}(x)$.
Assume that there is $y$ at level above $\alpha$ such that for each $t>\tau(x)$, $h_t$ is not constant on $S_{\ge y}$.
By definition of $\alpha$, this implies that $h_t(S_{\ge y})\cap U = \varnothing$ and
in particular that $h_t(y)\not\in U$ for each $t>\tau(x)$. But this contradicts the continuity of $h$ since $h_{\tau(x)}(y) = u\in U$.
Hence, $h_t$ is constant on $S_{\ge y}$ for at least one $t>\tau(x)$ and thus $\tau(y)>\tau(x)$.
This finishes the proof.
\endproof

Theorem \ref{thm:main-w1cpct} follows almost immediately.
\proof[Proof of Theorem \ref{thm:main-w1cpct}]
   Suppose that $W$ is a CW-complex and $f:M\to W$, $g:W\to M$ are such that $g\circ f$ is homotopic to $id_M$.
   If $f(S)$ meets uncountably many open cells, then $S$ contains an uncountable closed discrete
   subset (see for instance \cite[Theorem 1.5.19]{FritschPiccinini}) which contradicts its $\omega_1$-compactness.
   Hence $f(S)$ meets the interior of at most countably many cells and is therefore Lindel\"of.
   It follows that $g\circ f(S)$ is contained in a Lindel\"of hence metrizable submanifold of $M$.
   By Lemma \ref{lemgen:1} (i), $g\circ f(S)$ is constant on
   $S_{\ge x}$ if $x$ is above some fixed level $\alpha$ in $T$. Since $S_{\ge x}$ is uncountable for some $x$, 
   Theorem \ref{thm:main2} yields a contradiction.
\endproof
%%%%%%%%

\subsection{Proof of Theorem \ref{thm:main3}}

The proof is almost exactly the same as that of 
Theorem \ref{thm:main2} once we have an equivalent of Corollary \ref{cor:steprlocmet} in this context.
This is given by Corollary \ref{cor:omega_1} below.
Lemma \ref{lemma:omega_1} below plays the role of Lemma \ref{lemgen:1} (i). 
We first state the following easy fact.

\begin{lemma}\label{lemma:staclub}
   Let $S\subset\omega_1$ be stationary, endowed with the subspace topology.
   Then an at most countable family of club subsets of $S$ has a club intersection. 
\end{lemma}
\proof
   A direct proof is straightforward, but notice that
   since
   $S$ is an $\omega_1$-compact subspace of the tree $\omega_1$, the result is also a consequence of Lemma \ref{lemgen:1} (f)--(g).
\endproof

\begin{lemma}\label{lemma:omega_1}
  Let $S$ be a stationary subset of $\omega_1$ endowed with the subspace topology.
  If $Y$ is regular, $[\aleph_0,\aleph_1]$-compact and has $G_\delta$ points, then
  any continuous $f: S \to Y$ is eventually constant, that is, there is $\alpha\in\omega_1$ such that $f(\beta) = f(\alpha)$ for each $\beta\ge\alpha$,
  $\beta\in S$.
\end{lemma}

\proof 
  We start by showing that there is some $c\in Y$ such that $f^{-1}(\{c\})$ is club in $S$.
  If $f(S)$ is countable, then this is immediate.
  We thus assume that $f(S)$ is uncountable.
  Since $Y$ is $[\aleph_0,\aleph_1]$-compact, $f(S)$ has a condensation point $c\in Y$. 
  Since $Y$ is regular and has $G_\delta$ points,  we may choose 
  open sets $U_n\ni c$, $n\in\omega$, such that $\cap_{n\in\omega}U_n = \cap_{n\in\omega}\wb{U_n} =\{c\}$.
  Since $c$ is a condensation point, $f^{-1}(\wb{U_n})$ is club in $S$ for each $n$, hence
  $f^{-1}(\{c\})=f^{-1}(\cap_{n\in\omega}\wb{U_n}) = \cap_{n\in\omega}f^{-1}(\wb{U_n})$ is club in $S$ by Lemma \ref{lemma:staclub}.
  \\
  Now, since $f^{-1}(Y - U_n)$ is closed, it must be bounded, otherwise it intersects $f^{-1}(\{c\})$.
  It follows that $f^{-1}(Y-\{c\}) = \cup_{n\in\omega} f^{-1}(Y - U_n)$ is bounded in $S$, say by $\alpha$.
  Hence $f$ is constant on $S$ above $\alpha$.  
\endproof

\begin{cor}\label{cor:omega_1}
   Let $S$ be a stationary subset of $\omega_1$ endowed with the subspace topology.
   Let $Y$ be a regular space with $G_\delta$ points. Let $U\subset Y$ be open such that
   $\wb{U}$ is $[\aleph_0,\aleph_1]$-compact. Let $h:S\times[0,1]\to Y$ be continuous.
   Then there is $\alpha\in\omega_1$ such that for each $t\in[0,1]$ 
   either $h_t^{-1}(U) \cap [\alpha,\omega_1)\cap S = \varnothing$, or $h_t$ is constant on $[\alpha,\omega_1)\cap S$.
\end{cor}
\proof
   Fix a countable dense subset $\{t_n\,:\,n\in\omega\}\subset[0,1]$.
   By Lemmas \ref{lemma:staclub}--\ref{lemma:omega_1} we may fix $\alpha$ such that for each $n$
   either $[\alpha,\omega_1)\cap S \cap h_{t_n}^{-1}(\wb{U})=\varnothing$ or $h_{t_n}$ is constant above $\alpha$.
   The result follows by continuity.
\endproof

The proof of Theorem \ref{thm:main3} can now be done exactly along the same lines as that of Theorem \ref{thm:main2}, we thus only give a sketch.
Set $\tau = \sup\{t\,:\,h_t\text{ is eventually constant}\}$. 
Then $h_\tau$ is eventually constant and $\tau < 1$ since $h_1$ has uncountable image. 
Fix $\alpha$ minimal in $S$ such that $h_\tau$ is constant above $\alpha$, take an open $U$ containing $h_\tau(\alpha)$
such that $\wb{U}$ is $[\aleph_0,\aleph_1]$-compact. By Corollary \ref{cor:omega_1}
this contradicts the continuity of $h$.

%%%%%%%%%%%%%%%%%%%%%%%%%%%%%%%%%%%%%%%%%%%%%%%%%%%%%%%%

%%%%%%%%%%%%%%%%%%%%%%%%%%%%%%%%%%%%%%%%%%%%%%%%%%%%%%%%
%%%%%%%%%%%%%%%%%%%%%%%%%%%%%%%%%%%%%%%%%%%%%%%%%%%%%%%%

%\setcounter{section}{0}
%\renewcommand\thesection{\Alph{section}}

\section{Surfaces from trees}\label{sec:surfaces}

The goal of this section is to define a surface $M_T$ homotopy equivalent to the road space $R_T$ of a 
given rooted binary $\omega_1$-tree, and to check some additional properties involving the mapping class group.
We recall that we assume that our trees are Hausdorff.
(Everything in this section could be done starting from any rooted $\omega_1$-tree, but for simplicity we shall consider only 
binary ones.) 

\subsection{A surface homotopy equivalent to the road space of a tree}\label{sec:A1}

\begin{figure}[h]
  \begin{center}
  \epsfig{figure=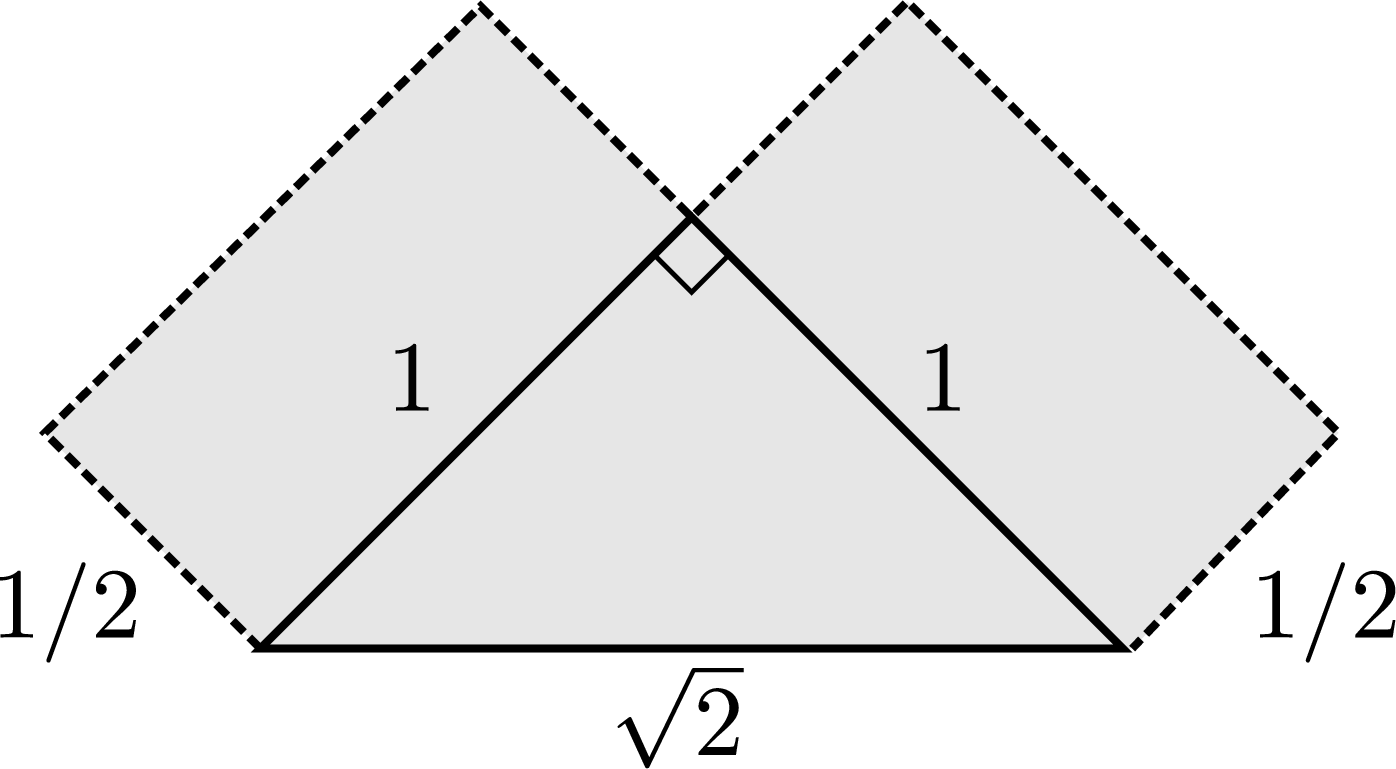,height=3.5cm}
  \caption{The `pieces' $M_x$.}
  \label{fig1}
  \end{center}
\end{figure}

\begin{figure}[h]
  \begin{center}
    \epsfig{figure=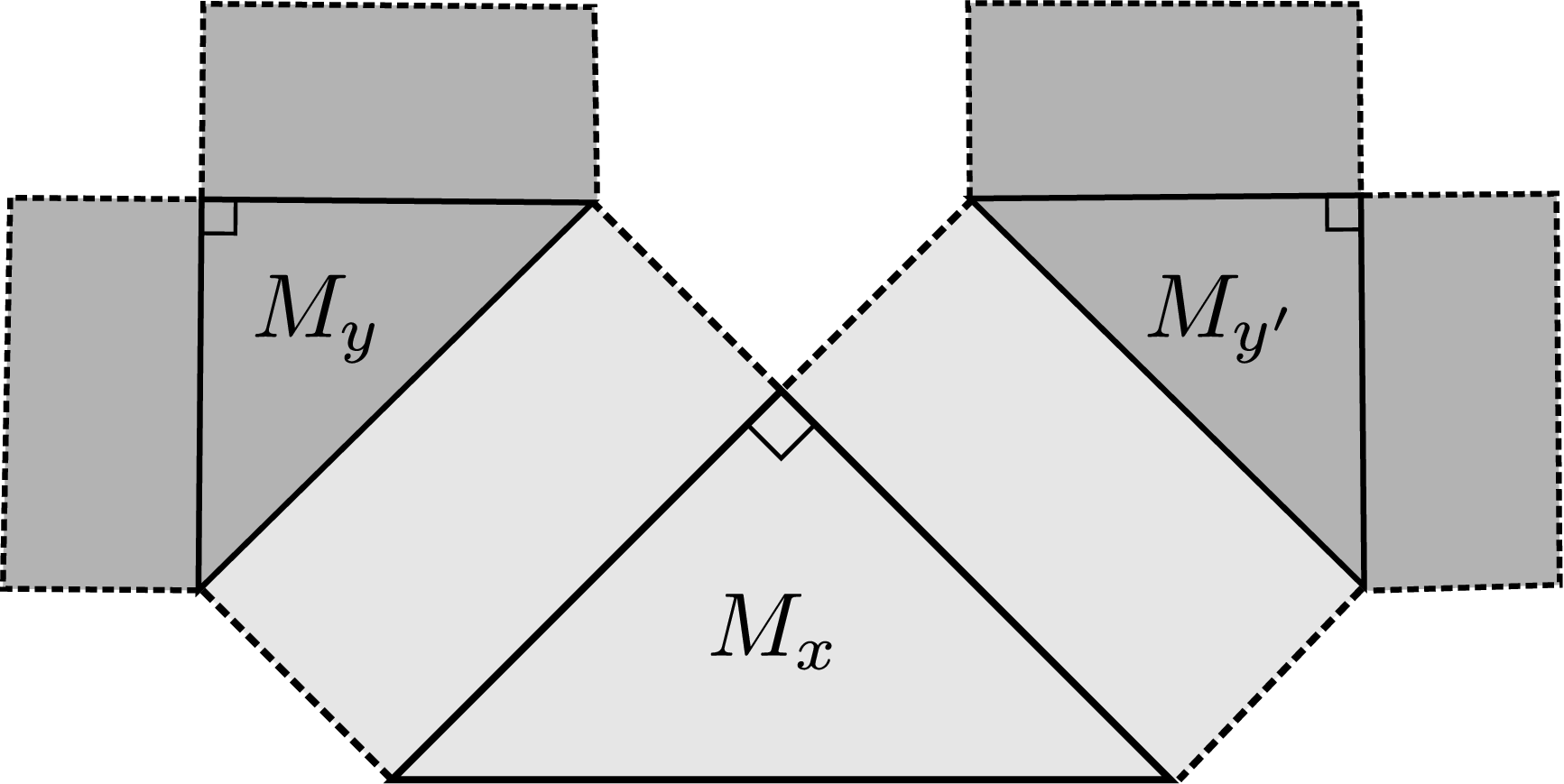, height=3.8cm}
    \caption{Gluing successor levels}
    \label{fig2}
  \end{center}
\end{figure}

This construction was found and explained to us by Peter Nyikos in 2006.
We fix a binary $\omega_1$-tree $T$. Let us define the manifold $M_T$. 
As a set, $M_T=\cup_{x\in T}M_x$, where each piece $M_x$ is the union of a right isoceles triangle
with sides $1,1,\sqrt{2}$ and two rectangles with sides $1,\frac{1}{2}$ glued together as in Figure \ref{fig1}
(dashed boundaries do not belong to $M_x$).
The topology in the interior of $M_x$ is the usual topology. We shall need an arbitrary total order $\vartriangleleft$ on $T$.
If $y,y'$ are the two immediate $\le$-successors of $x$ in $T$ and $y\vartriangleleft y'$, 
we glue $M_y$ and $M_{y'}$ on the top of $M_x$ (after multiplying them by $1/\sqrt{2}$) as in Figure \ref{fig2}
below
(we put the $\vartriangleleft$-bigger on the right), which gives a topology for points
in the bottom boundary of $M_x$ for $x$ at successor levels.

\begin{figure}[h]
  \begin{center}
  \epsfig{figure=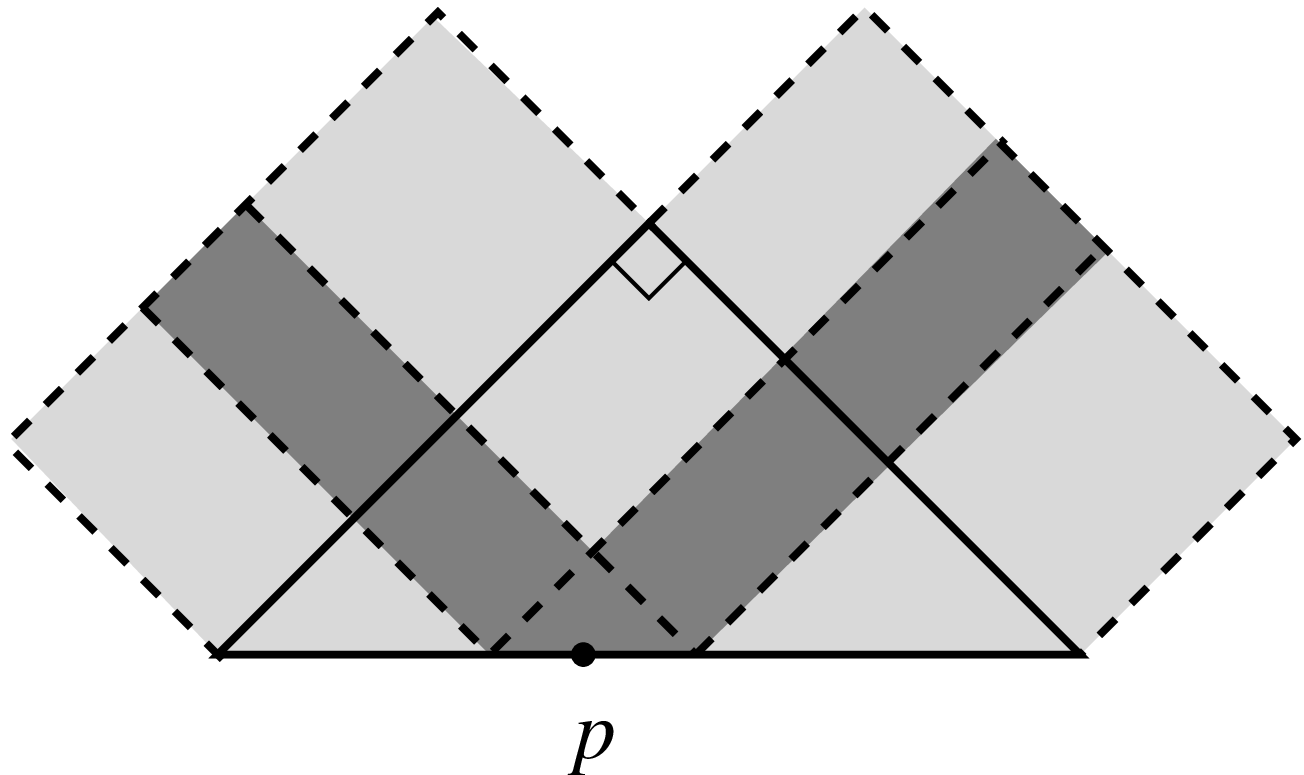,height=3.5cm}\hskip 1cm
  \epsfig{figure=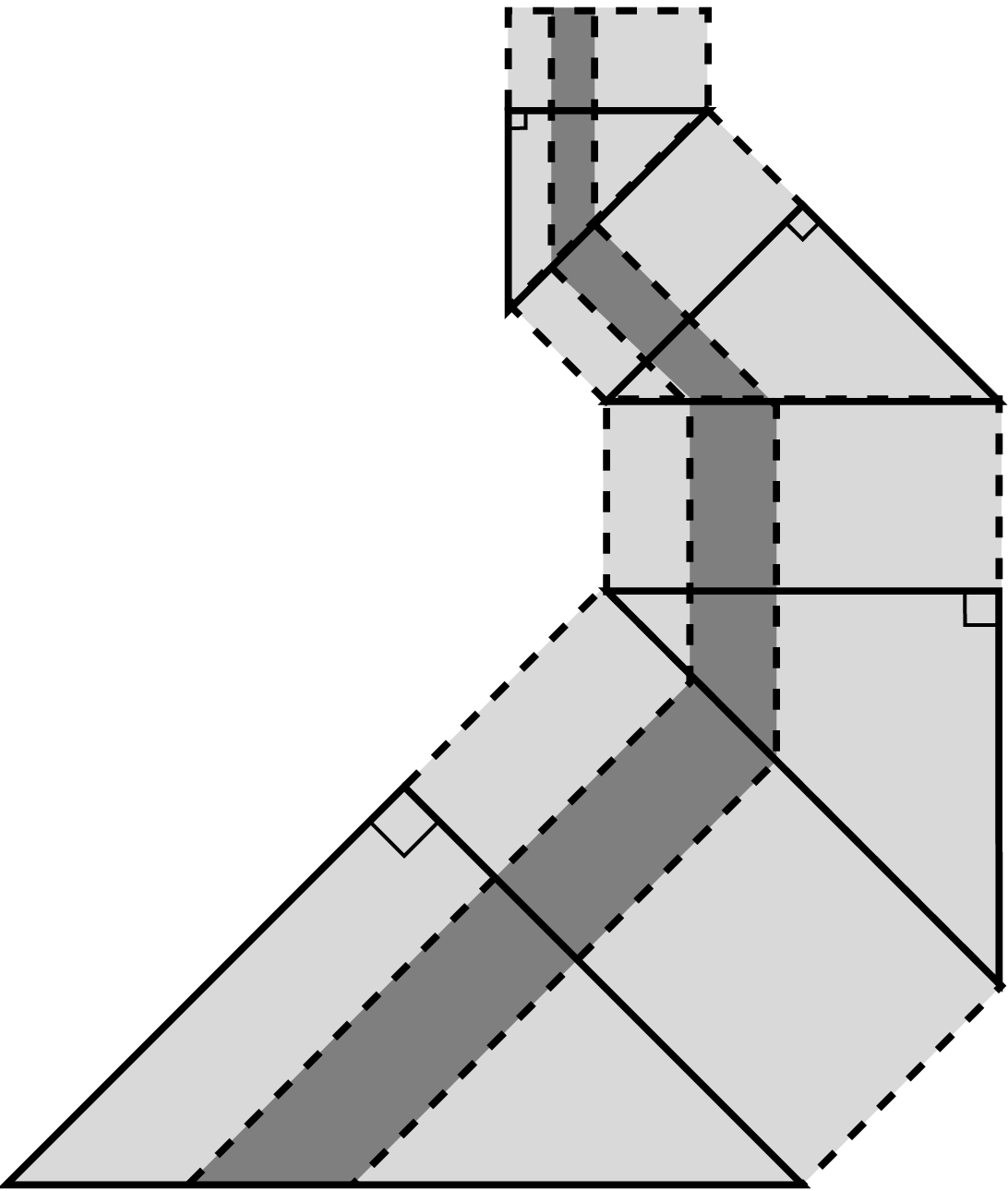,height=3.8cm}
  \caption{Neighborhoods at limit levels}
  \label{fig3}
  \end{center}
\end{figure}

If now $p\in M$ belongs to the bottom boundary of $M_x$ for $x$ at a limit level $\alpha$ of $T$, a
neighborhood of $p$ is given by a half ball in $M_x$ of diameter $\epsilon$ centered at $p$, together with
a `zigzag' running between some level $\gamma<\alpha$ and $\alpha$ which follows $p$ in $M_T$,
that is, the union of open strips of width $\epsilon$ in $M_y$ for all $y<x$ 
in levels $\ge\gamma$ (see Figure \ref{fig3}), 
choosing the strip that points eastwards or westwards according to $\vartriangleleft$. Then,
to make it open,
remove the part that lies in the bottom boundary of $M_y$, for $y\in\text{\rm Lev}_\gamma(T)$.

\begin{figure}[h]
  \begin{center}
     \epsfig{figure=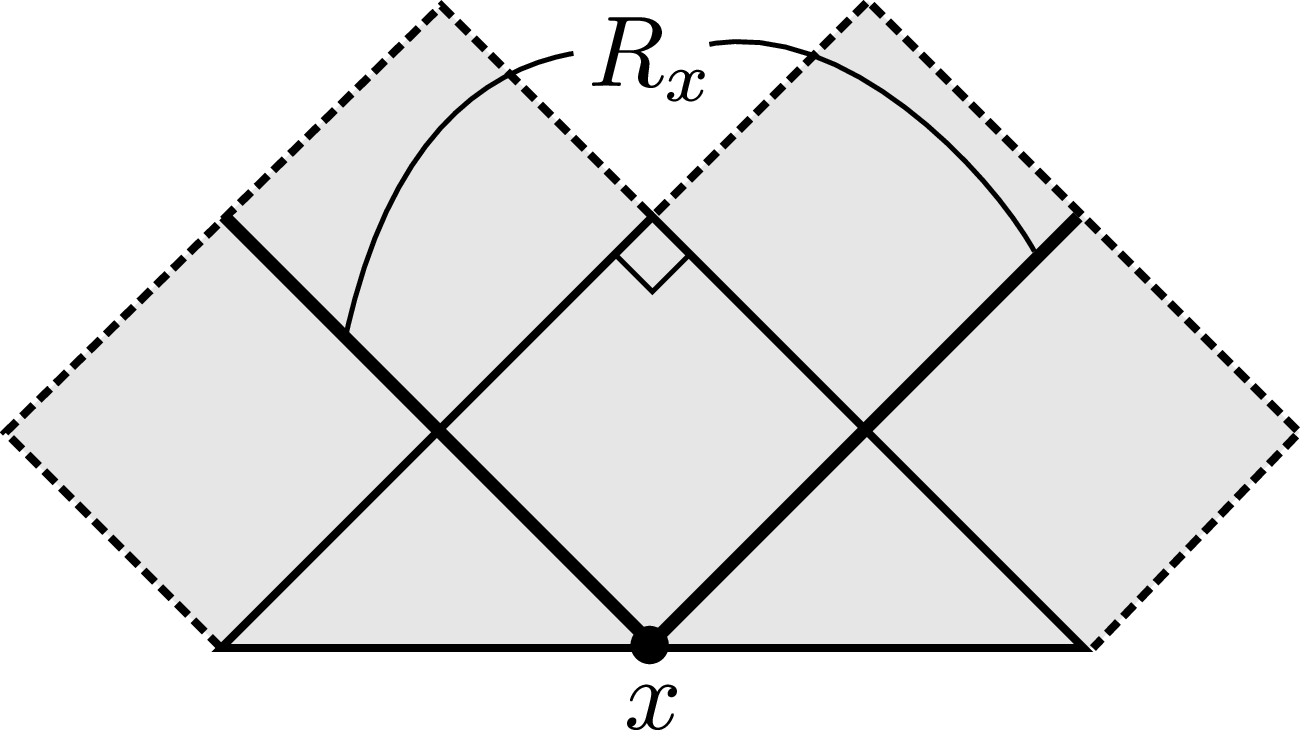,height=3.2cm}
     \caption{$T$ and $R_T$ in $M_T$.}
     \label{fig4}
  \end{center}
\end{figure}

This gives $M_T$ a structure of Type I surface without boundary 
with canonical cover given by $M_\alpha=\cup_{\beta<\alpha}\cup_{x\in\text{\em Lev}_\beta(T)}M_x$.
(Strictly speaking, there is a remaining boundary in $M_0$, $0$ the root of $T$, but we remove it
by adding an extra open strip below it.)
Since a countable increasing union of $2$-discs is homeomorphic to a $2$-disc 
(the famous M. Brown theorem, here in dimension $2$, see for instance
\cite[Section 3.1]{GauldBook}), by induction $M_\alpha$ is homeomorphic to a $2$-disc, and $M_T$ is a surface.
$T$ and $R_T$ embed in $M_T$ as seen in Figure \ref{fig4} (where $R_x=R\cap M_x$). It is clear that all homotopy groups of
$M_T$ are trivial, because the image of a compact set by a continuous map will be contained in some $M_\alpha$.
Notice that $T$ is closed in $M_T$ but $R_T$ is not: if $x$ is at a limit level in $T$,
then at least one of the left or right half part of the bottom boundary of $M_x$ is in the closure of the segments 
in $R_T$ emanating
from the points of $T$ below $x$.

\begin{thm}\label{propcontractible} (P. Nyikos)
   The inclusion $R_T\subset M_T$ is a homotopy equivalence.
\end{thm}

This theorem is a direct consequence of the following lemma (see Exercise 4 p. 18 in \cite{Hatcher}): 

\begin{lemma}
  \label{lemma2}
  There is a deformation retraction in the weak sense from $M_T$ to $R_T$, i.e. there is a homotopy
  $h:M_T\times[0,1]\to M_T$ with $h_0=\text{\rm id}$, $h_1(M_T)=R_T$ and $h_t(R_T)=R_T$ for all $t\in[0,1]$.
\end{lemma}

The geometric features of $M_T$ will be used in the proof.
It is interesting to note that
it is not possible to obtain a (real) deformation retraction 
(i.e. $h_t\upharpoonright R_T =\text{\rm id}_{R_T}$ for each $t$), in fact there is even no retraction $h:M_T\to R_T$ (i.e. with 
$h\upharpoonright R_T =\text{\rm id}_{R_T}$), since $R_T$ is not closed in $M_T$.
 
\begin{figure}[h]
  \begin{center}
    \epsfig{figure=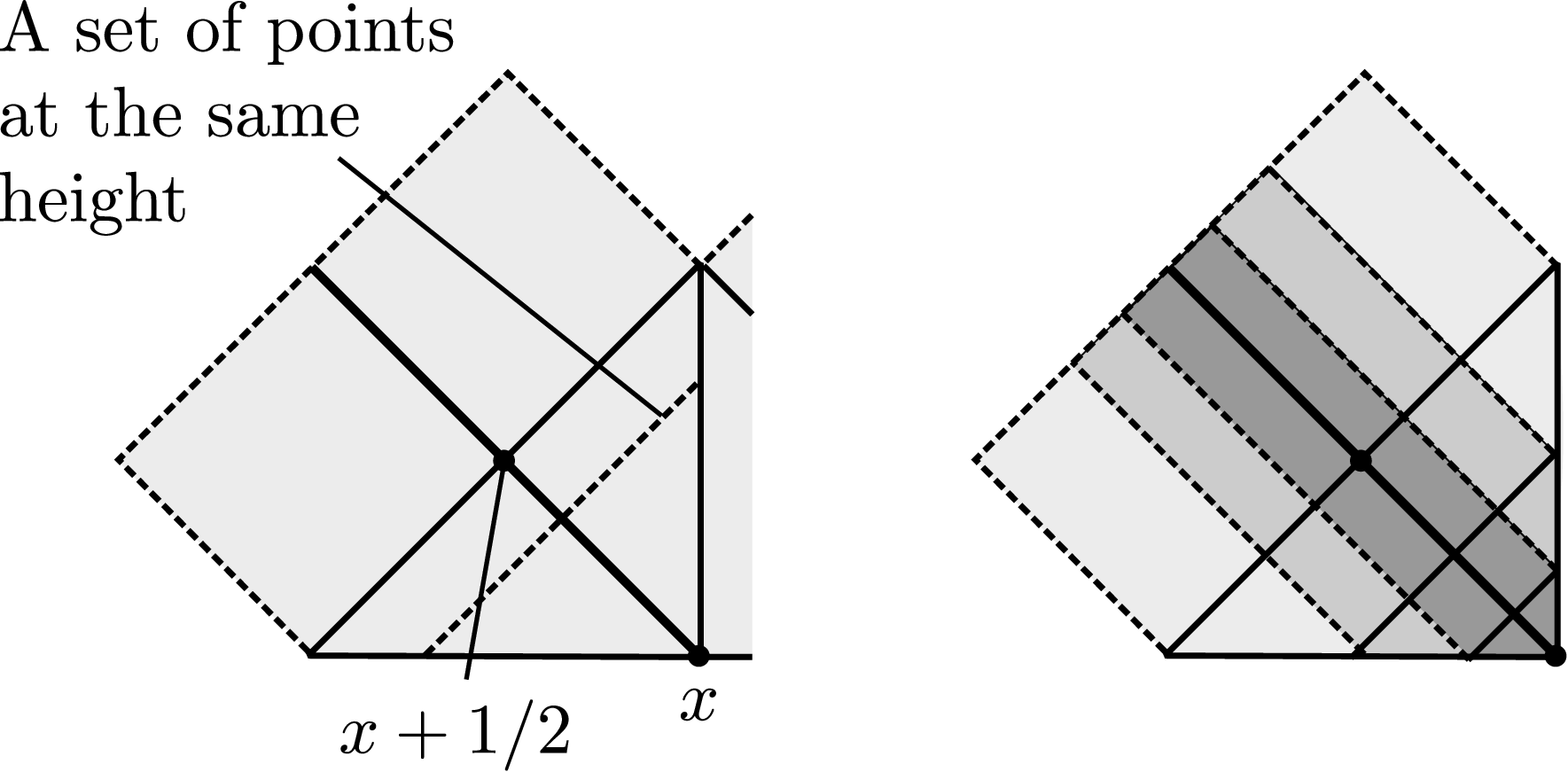,height=4cm}
    \caption{The homotopy globally seen.}
    \label{fig5}
  \end{center}
\end{figure}

\proof[Proof of Lemma \ref{lemma2}.]
The homotopy we shall describe is the same in every piece $M_x$ of $M_T$, and acts symmetrically on $M_x$, 
we thus only need to explain what happens in its left part (see Figure \ref{fig5}).
Let us first describe how the homotopy works globally. The points in $R_x$ will stay in it, those on the
segment between $x$ and $x+1/2$ (see Figure \ref{fig5}, left, for the definition) 
will be all contracted to the point $x$ while the segment between 
$x+1/2$ and
the next level will eventually cover all of $R_x$. Meanwhile, the south-west and north-east `sides' will
shrink at unit speed to $R_x$, and the triangle will be completely contracted to $x$ (see Figure \ref{fig5},
on the right).

\begin{figure}[h] 
  \begin{center}
    \epsfig{figure=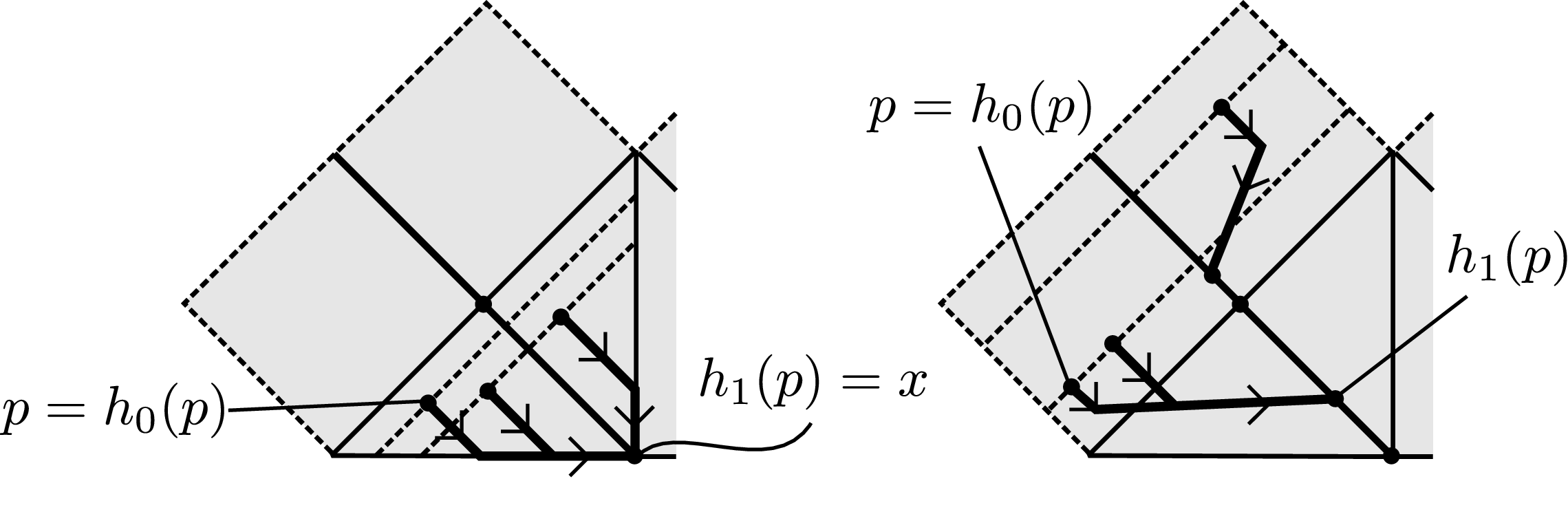,height=4cm}
    \caption{Trajectories under $h_t$.}\label{fig6}
  \end{center}
\end{figure}

Let us now be more specific and describe the trajectories of the points under the homotopy.
First, the triangle.
The points at height $x+1/2$ start moving at once, while the points in the triangle start moving when
the points at height $x+1/2$ `reach them'. Then, they go at 
speed $\frac{1}{2}$ parallel to $R_x$ until they hit
the bottom or the vertical line above $x$, and then they head directly toward $x$ at speed 
$\frac{\sqrt{2}}{2}$ (see Figure
\ref{fig6}, on the left). Hence, $h_1(p)=x$ for all points in the triangle.
Points in the rectangle have one velocity component parallel to $R_x$ which is proportional
to the distance to the next level (that is, between $0$ and $\frac{1}{2}$), 
the other component is at speed $\frac{1}{2}$ and goes perpendicular to $R_x$.
The first component starts immediately, but the second one only when the south-west or the north-east side
reaches it (see Figure \ref{fig6}, on the right). 
The homotopy $h_t$ described above is continuous with
$h_0=id$ and $h_1(M_T)=R_T$. 
\endproof

Theorem \ref{thm:main-contractible} follows immediately:
\proof[Proof of Theorem \ref{thm:main-contractible}]
Let $T$ be a binary $\R$-special $\omega_1$-tree, and set $M=M_T$. Then, by Theorems
\ref{propcontractible} and \ref{thm:main1}, $M$ is a Type I contractible surface. 
\endproof

\begin{figure}[h] 
  \begin{center}
    \epsfig{figure=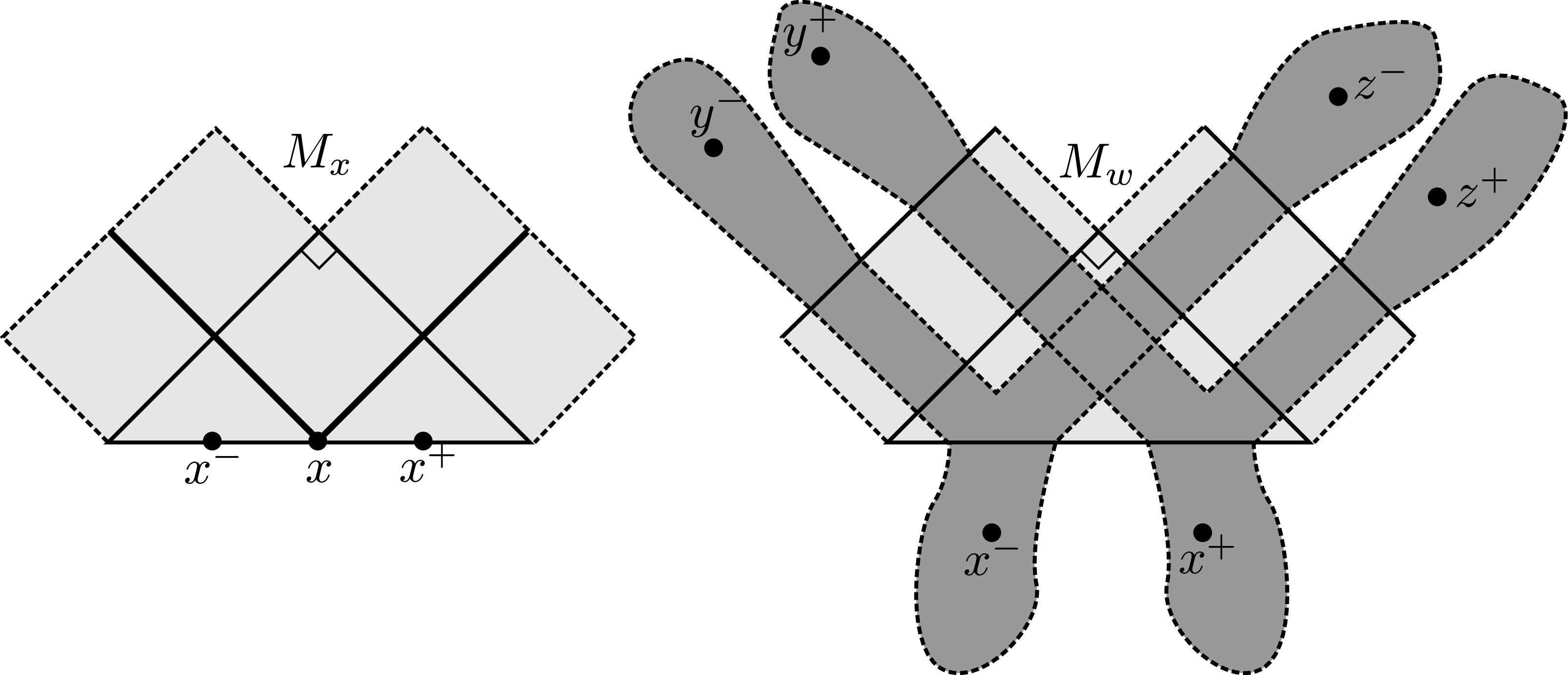,width=.8\textwidth}
    \caption{Lemma \ref{lemma:MTnotnormal}.}
    \label{fig:LemmaMT}
  \end{center}
\end{figure}

The tree $T$ embeds in $M_T$ as a closed subset.
If $T$ is $\Q$-special (that is, there is a strictly increasing function $T\to\Q$),
then $T$ is normal in some models of set theory, for instance under {\bf MA + $\neg$CH},
see e.g. \cite{Nyikos:trees}.
We might hope that normality `lifts up' to $M_T$,
but this is never the case if $T$ is Aronszajn.

\begin{lemma}\label{lemma:MTnotnormal}
   If $T$ is an Aronszajn tree, then $M_T$ is not normal.
\end{lemma}
The proof is reminiscent of Lemma 4.9 in \cite{Nyikos:1992}.
\proof
       For each $x\in T$, let $x^-,x^+\in M_T$ be defined as 
       on the lefthandside of Figure \ref{fig:LemmaMT}. Define $E^\pm=\{x^\pm\,:\,x\in T\}$,
       then $E^+,E^-$ are closed subsets of $M_T$ and have empty intersection.
       Let $U^+\supset E^+,U^-\supset E^-$ be open.
       For each $x\in T$, choose $U^+_x,U^-_x$ to be basic connected open subsets of $U^+,U^-$ containing
       respectively $x^+$ and $x^-$. If $x$ is at limit level $\alpha$,
       we may assume that $U^+_x,U^-_x$ are union of strips as in Figure \ref{fig3}, hence
       there is $\gamma(x)<\alpha$ 
       such that $y^+\in U^+_x$ and $y^-\in U^-_x$, where $y=x\upharpoonright\gamma(x)$.
       By the pressing down lemma for $\omega_1$-trees (see e.g. \cite[p. 154]{Hart-Souslin}),
       there is $\gamma\in\omega_1$ and
       $S\subset T$ meeting a stationary set of levels such that $\gamma(x) = \gamma$ for each $x\in S$.
       Since $\lev_\gamma(T)$ is countable, there is some $y\in T$ at level $\gamma$ such that
       $y^+\in U^+_x$ and $y^-\in U^-_x$ for $x$ in an uncountable set $S'\subset T$.
       It follows that both $V^+ = \cup_{x\in S'}U^+_x$ and $V^- = \cup_{x\in S'}U^-_x$
       are unbounded connected subsets of $U^+$ and $U^-$, respectively.
       Since $T$ is Aronzsajn and $S'$ is uncountable, there are $x,y,z$ in $S'$
       with $x<y,x<z$. Since $T$ is Hausdorff there is $w$ at maximal level
       such that $w<y,w<z$.  
       By connectedness, either $V^+$ intersects $V^-$ below the level of $w$,
       or they intersect in $M_w$, as shown on the righthandside of Figure \ref{fig:LemmaMT}.
       This shows that $U^+\cap U^-\not=\varnothing$, hence $M_T$ is not normal.
\endproof

%%%%%%%%%%%%%%%%%%%%%%%%%%%%

\subsection{Isotopies and the mapping class group}\label{sec:mapping}

Our goal is to prove the following theorem.

\begin{thm}\label{thm:mapping}
   \ \\
   (1) The mapping class group of the Pr\"ufer surface $\mathsf{P}$ (which is contractible) has cardinality $\ge\mathfrak{c}$.\\
   (2) There is a contractible Type I surface with mapping class group of cardinality $\ge\mathfrak{c}$.
\end{thm}

Recall that two homeomorphisms are {\em isotopic} if there is a homotopy $h_t$ between them which is an homeomorphism
for each $t\in[0,1]$.
The {\em mapping class group} $\text{\rm MG}(M)$ of a topological space $X$ is the group of homeomorphisms of $X$ quotiented by 
the isotopy classes.
Theorem \ref{thm:mapping} shows a sharp contrast between metrizable and non-metrizable surfaces:
two homeomorphisms of a metrizable surface are homotopic if and only if they are isotopic 
(except for some trivial exceptions). To our knowledge this is due to D.B.A Epstein \cite{Epstein:1966}.

We shall prove (2) in detail, and then explain why (1) holds as well.
We first recall some definitions from \cite{Nyikos:1984}.
Given a non-metrizable Type I manifold $M$ with canonical cover $\langle M_\alpha\,:\,\alpha\in\omega_1\rangle$, 
the tree of non-metrizable-component 
boundaries $\Upsilon(M)$ is defined as  
the set of topological boundaries $\partial C$ such that $C$ is a non-metrizable connected component of $M - M_\alpha$ 
for some $\alpha\in\omega_1$, with the following order:
$\sigma\ge\tau$ if $\sigma$ is a subset of a component whose boundary is $\tau$.
(Of course, $\Upsilon(M)$ depends on the choice of the canonical cover, but this is not inconsequential for our 
considerations.)
Recall that an $\omega_1$-tree $T$ is {\em normalized} if $T_{\ge x}$ is uncountable for each $x\in T$.
Since we only take non-metrizable components, $\Upsilon(M)$ is normalized,
and
if one starts with a normalized $T$, then 
$\Upsilon(M_T)=T$ for the surface $M_T$ constructed in Section \ref{sec:A1}.

Given a tree $T$, we say that a map $T\to T$ is an {\em automorphism} of $T$ if it is an order preserving bijection.
A map $f:M\to M$ of a Type I manifold $M$ 
(with some given canonical cover)
is a $\Upsilon$-automorphism if it is a homeomorphism of $M$ such that for all 
$\sigma\in\Upsilon(M)$, 
$f(\sigma)=\tau$ for some $\tau\in\Upsilon(M)$, and if the induced map $\wh{f}:\Upsilon(X)\to\Upsilon(X)$ is an
automorphism. Being a $\Upsilon$-automorphism is a strong assumption, however it is not difficult to 
see that any homeomorphism of a Type I manifold induces an $\Upsilon$-automorphism of $\Upsilon(M)$ for a
well choosen canonical cover.
 
\begin{lemma} \label{nonisotopic} 
  Let $M=\cup_{\alpha\in\omega_1}M_\alpha$ be a Type I manifold.
  If $f,g:M\to M$ are $\Upsilon$-automorphisms such that 
  $\wh{f}\not=\wh{g}$, then 
  $f$ and $g$ are non isotopic.
\end{lemma}
We shall not give the shortest proof of this lemma but rather an argument which can be easily adapted for the Pr\"ufer surface.

\iffalse
\proof 
  Since $X$ is of Type I, $X -  X_\alpha$ cannot have uncountably many connected components, otherwise
  $\wb{X_{\alpha+1}}$ would not be Lindel\"of. Hence $\Upsilon(X)$ is an $\omega_1$-tree.
  Notice that $\Upsilon(X)$ is also well pruned since a non-metrizable component of $X$ 
  must be unbounded in $X$.
  We can assume that $g=id$.
  If $h:X\to X$ is a continuous map, then $C(h)=\{\alpha\in\omega_1\,:\,h(X_\alpha)\subset X_{\alpha}\}$ is club.
  If $h_t$ is an isotopy between $f$ and the identity, 
  then $C=\bigcap_{t\in[0,1]\cap\Q}C(h_t)\cap C(h_t^{-1})$ is also club. For $\alpha\in C$,
  we have $h_t(\wb{X_\alpha} -  X_\alpha)=\wb{X_\alpha} -  X_\alpha$ for all $t$. It follows that 
  for all $\alpha\ge\min C$, the components of $X -  X_\alpha$ do not move under $h_t$. Since the tree is normalized, 
  no component can move at all and $\wh{f}=id$. 
\endproof
\fi

\proof
  Let 
  $h_t$ be a homotopy such that $h_0=f$, $h_1=g$. Let $\sigma\in \Upsilon(M)$ be such that 
  $\wh{f}(\sigma)\not= \wh{g}(\sigma)$. We will show that $h_t$ cannot be an isotopy.
  Since
  $\wh{f},\wh{g}$ are automorphisms, $\wh{f}(\tau)\not= \wh{g}(\tau)$ for all $\tau\ge\sigma$
  in $\Upsilon(M)$, and
  both $\wh{f}(\sigma)$ and $\wh{g}(\sigma)$ belong to the same level $\gamma$ of $\Upsilon(M)$. By
  definition of $\Upsilon(M)$,
  there is some $\alpha$ such that $\wh{f}(\tau)$ and $\wh{g}(\tau)$ are in 
  different connected components of 
  $M -  M_\alpha$ for all $\tau>\sigma$. Thus, for $\tau\ge\sigma$ and all
  $y\in f(\tau)\subset M$,
  there is an open (connected) interval
  $I_y\subset [0,1]$ such that $h_t(y)\in M_\alpha$ for $t\in I_y$. 
  Now, since $\Upsilon(M)$ is normalized, for all $\beta\in\omega_1$ with $\beta>\gamma$, 
  we may choose $\tau_\beta$ in the $\beta$-level of $\Upsilon(M)$ and 
  $y_\beta$ in $f(\tau_\beta)$.
  Since $[0,1]$ is second countable, there
  is an interval $J$ contained in $I_{y_\beta}$ 
  for uncoutably many $\beta$.
  Thus, for $t\in J$, $h_t$ sends an unbounded and therefore non-metrizable 
  subset of $M$ inside $M_\alpha$,
  which is metrizable. Hence, $h_t$ cannot be a homeomorphism, and $h_t$ is not an isotopy. 
\endproof

\proof[Proof of Theorem \ref{thm:mapping}.]
We first prove (2).
Let $T$ be any rooted binary $\omega_1$-tree and $M_T\supset T$ 
be as in Section \ref{sec:A1}.
Take now copies $T_i\subset M_i$ of $T,M_T$ for each $i\in\Z$. Build the manifold
$M^\omega$ by gluing all these $M_i$ on a open disc as in Figure \ref{fig7} below.
Let $T^\omega$ be $\cup_{i\in\Z}T_i$ together with a new root $x$ as shown in Figure \ref{fig7}. 
For $m\in\Z$ let $\varphi_m:M^\omega\to M^\omega$ be the homeomorphism acting 
by `rotating' $M_i$ to $M_{i+m}$. It is clear that $\varphi_m$ is a $\Upsilon$-automorphism
of $M^\omega$.
\begin{figure}[h]
  \begin{center}
    \epsfig{figure=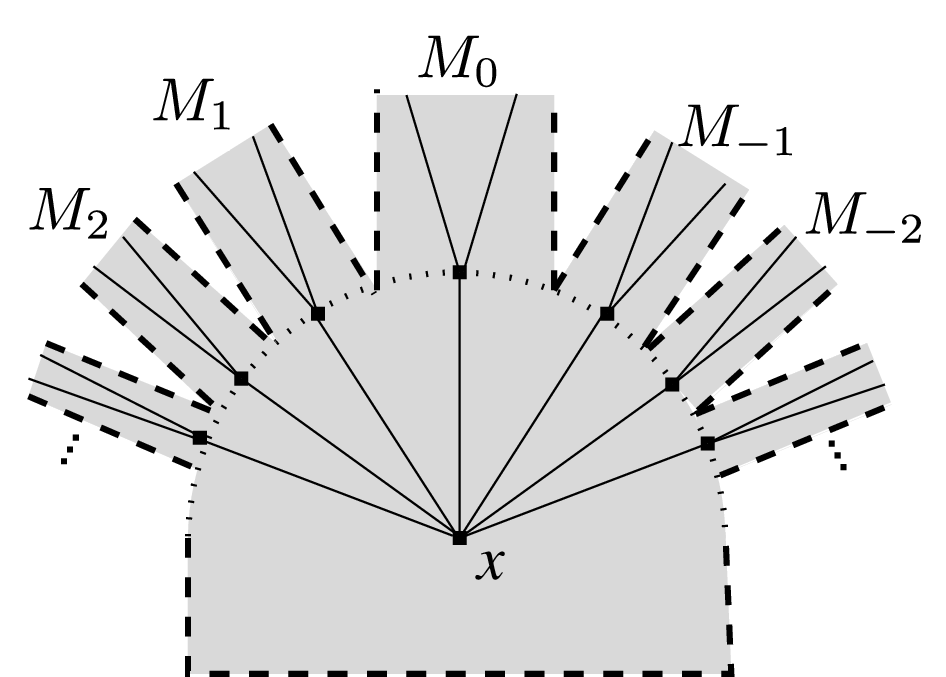,height=4cm}
    \caption{$M^\omega$.}
    \label{fig7}
  \end{center}
\end{figure}
We now define $M=M^{\omega,\omega}$ with copies $M_j^\omega$ of $M^\omega$ for $j\in\Z$ in the same
way. 
Again, $\Upsilon(M)=T^{\omega,\omega}$, where $T^{\omega,\omega}$ is the $\omega_1$-tree obtained from
the $T^\omega_i$ by adding a common root. 
For $s\in\Z^\Z$, let $\phi_s:M\to M$ be given
by $\phi_s\upharpoonright M^\omega_j =\varphi_{s(j)}:M^\omega_j\to M^\omega_j$ and $\phi_s(u)=u$ if
$u\not\in\cup_{j\in\Z}M^\omega_j$. It is clear that $\phi_s$ is an $\Upsilon$-automorphism.
By Lemma \ref{nonisotopic},
$s\mapsto [\phi_s]\in\text{\rm MG}(M)$ is one--to--one, hence $\text{\rm MG}(M)$ has cardinality at least $\mathfrak{c}$.
If $T$ is a rooted binary $\R$-special tree, then 
$M$ is contractible.
\\
We now briefly explain how to prove (1). 
We use the notations of Section \ref{sec:prufer} for everything related to the Pr\"ufer surface $\mathsf{P}$.
Given $a\in\R$,  let $f_a:\mathsf{P}\to\mathsf{P}$ send all of $\R_c$ to $\R_{c+a}$ pointwise,
and be defined by $\langle x,y\rangle\mapsto\langle x,y+a\rangle$ in $H_0$.
Then each $f_a$ is a homeomorphism.
For each $c\in\R$ take the $0$-point $0_c\in\R_c$.
The proof of Lemma \ref{nonisotopic} can be easily adapted to show that if $h$ is an homotopy between 
$f_a$ and $f_b$ (with $a\not= b$),
then there is some $t$ such that 
$h_t$ sends $0_c$ inside $H_0$ for an uncountable subset $A\subset\R$ of $c$. 
Recall that $\{0_c\,:\,c\in A\}$ is closed discrete for any $A$.
Since $H_0$ is hereditarily separable,
the image of these points cannot be closed discrete, so $h_t$ is not an homeomorphism and $h$ is not an isotopy.
This shows that the map $\R\to\text{\rm MG}(\mathsf{P})$ defined by $a\mapsto [f_a]$ is one-to-one.
\endproof

The construction described in the proof of Point (2) but using the first octant in the square of the longray instead
of $M_T$ yields an $\omega$-bounded surface with mapping class group of cardinality $\ge\mathfrak{c}$,
see \cite[Thm 5.31]{GauldBook}.


\begin{thebibliography}{10}

\bibitem{AlexandroffUrysohn}
P.~Alexandroff and P.~Urysohn.
\newblock {\em {M\'emoire} sur les espaces topologiques compacts}.
\newblock Number~14 in Proceedings of the section of mathematical sciences.
  Koninklijke Nederlandse Akademie van Wetenschappen te Amsterdam, 1929.

\bibitem{meszigueshom}
M.~Baillif.
\newblock The homotopy classes of continuous maps between some non-metrizable
  manifolds.
\newblock {\em Topology Appl.}, 148(1--3):39--53, 2005.

\bibitem{mesziguessurf}
M.~Baillif.
\newblock Homotopy in non-metrizable $\omega$-bounded surfaces.
\newblock preprint arXiv:math/0603515v2, 2006.

\bibitem{MesziguesBridges}
M.~Baillif.
\newblock Some steps on short bridges: non-metrizable surfaces and
  {CW}-complexes.
\newblock {\em Journal of Homotopy and Related Structures}, 7(2):153--164,
  2012.

\bibitem{meszigues-contract-trees}
M.~Baillif.
\newblock (non-)contractible road spaces of trees, 2019.
\newblock arXiv preprint 1908.09901.

\bibitem{mesziguesnarrow}
M.~Baillif.
\newblock Relative (functionally) {Type I} spaces and narrow subspaces.
\newblock {I}n preparation, 2022.

\bibitem{BaillifDeoGauld}
M.~Baillif, S.~Deo, and D.~Gauld.
\newblock The mapping class group of powers of the long ray and other
  non-metrisable surfaces.
\newblock {\em Topology Appl.}, 157(8):1314--1324, 2010.

\bibitem{Balogh:1989}
Z.~Balogh.
\newblock On compact {H}ausdorff spaces of countable tightness.
\newblock {\em Proc. Amer. Math. Soc.}, 105:755--764, 1989.

\bibitem{BaloghRudin:1992}
Z.~Balogh and M.E. Rudin.
\newblock Monotone normality.
\newblock {\em Topology Appl.}, 47(2):115--127, 1992.

\bibitem{CalabiRosenlicht}
E.~Calabi and M.~Rosenlicht.
\newblock Complex analytic manifolds without countable base.
\newblock {\em Proc. Amer. Math. Soc.}, 4:335--340, 1953.

\bibitem{DeoGauld:2007}
S.~Deo and D.~Gauld.
\newblock Eventually constant spaces and nonmetrizable homology spheres.
\newblock {\em J. Ind. Math. Soc.}, Special centenary Volume
  (1907--2007):167--175, 2007.

\bibitem{DevlinShelah}
K.~J. Devlin and S.~Shelah.
\newblock Souslin properties and tree topologies.
\newblock {\em Proc. London Math. Soc.}, 39:237--252, 1979.

\bibitem{EisworthNyikos}
T.~Eisworth and P.~Nyikos.
\newblock Antidiamonds principles and topological applications.
\newblock {\em Trans. Amer. Math. Soc.}, 361:5695--5719, 2009.

\bibitem{Engelking}
R.~Engelking.
\newblock {\em General topology}.
\newblock Heldermann, Berlin, 1989.
\newblock Revised and completed edition.

\bibitem{Epstein:1966}
D.B.A Epstein.
\newblock Curves on 2-manifolds and isotopies.
\newblock {\em Acta Math.}, 115:83--107, 1966.

\bibitem{FritschPiccinini}
R.~Fritsch and R.A. Piccinini.
\newblock {\em Cellular structures in topology}, volume~19 of {\em Cambridge
  studies in advanced mathematics}.
\newblock Cambridge University Press, 1990.

\bibitem{GabardWouuuh}
A.~Gabard.
\newblock A separable manifold failing to have the homotopy type of a
  {CW}-complex.
\newblock {\em Arch. Math.}, 90:267--274, 2008.

\bibitem{GauldBook}
D.~Gauld.
\newblock {\em Non-metrisable manifolds}.
\newblock Springer-Verlag, New York-Berlin, 2014.

\bibitem{Greenwood:2002}
S.~Greenwood.
\newblock Constructing type {I} nonmetrisable manifolds with given
  {$\Upsilon$}-trees.
\newblock {\em Topology Appl.}, 123(1):91--101, 2002.
\newblock Proceedings of the Janos Bolyai Mathematical Society 8th
  International Topology Conference (Gyula, 1998).

\bibitem{GreenwoodNyikos:2005}
S.~Greenwood and P.~Nyikos.
\newblock {$\omega_1$}-compactness in {Type I} manifolds.
\newblock {\em Topology Appl.}, 48(1--3):165--171, 2005.
\newblock Proceedings of the Janos Bolyai Mathematical Society 8th
  International Topology Conference (Gyula, 1998).

\bibitem{Hart-Souslin}
K.P. Hart.
\newblock More remarks on {S}ouslin properties and tree topologies.
\newblock {\em Top. App.}, 15:151--158, 1983.

\bibitem{Hatcher}
A.~Hatcher.
\newblock {\em Algebraic topology}.
\newblock Cambridge University Press, Cambridge, 2002.
\newblock (A free electronic version is available).

\bibitem{IsmailNyikos}
M.~Ismail and P.~Nyikos.
\newblock On spaces in which countably compact sets are closed, and hereditary
  properties.
\newblock {\em Topology Appl.}, 11(3):281--292, 1980.

\bibitem{kunen}
K.~Kunen.
\newblock {\em Set theory, an introduction to independance proofs (reprint)}.
\newblock Number 102 in Studies in Logic and the Foundations of Mathematics.
  North-Holland, Amsterdam, 1983.

\bibitem{MardaniThesis}
A.~Mardani.
\newblock Topics in the general topology of non-metric manifolds.
\newblock The University of Auckland. Ph.D Thesis, 2014.

\bibitem{MilnorCW}
J.~Milnor.
\newblock On spaces having the homotopy type of a {CW}-complex.
\newblock {\em Trans. Am. Math. Soc.}, 90:272--280, 1959.

\bibitem{Nyikos:1984}
P.~Nyikos.
\newblock The theory of nonmetrizable manifolds.
\newblock In {\em Handbook of {S}et-{T}heoretic {T}opology (Kenneth Kunen and
  Jerry~E. Vaughan, eds.)}, pages 633--684. North-Holland, Amsterdam, 1984.

\bibitem{Nyikos:1992}
P.~Nyikos.
\newblock Various smoothings of the long line and their tangent bundles.
\newblock {\em Adv. Math.}, 93(2):129--213, 1992.

\bibitem{Nyikos:trees}
P.~Nyikos.
\newblock Various topologies on trees.
\newblock In P.R. Misra and M.~Rajagopalan, editors, {\em Proceedings of the
  Tennessee Topology Conference}, pages 167--198. World Scientific, 1997.

\bibitem{NyikosZdomskyy}
P.~Nyikos and L.~Zdomskyy.
\newblock Locally compact, {$\omega_1$-compact} spaces, 2017.
\newblock ArXiv preprint 1712.03906, submitted.

\bibitem{Schlindwein:2003}
C.~Schlindwein.
\newblock How special is your {A}ronszajn tree?
\newblock eprint arXiv:math/0312445, 2003.

\bibitem{CEIT}
L.~A. Steen and J.~A.~Seebach Jr.
\newblock {\em Counterexamples in topology}.
\newblock Springer Verlag, New York, 1978.

\bibitem{StepransTrees}
J.~{Stepr\=ans}.
\newblock Trees and continuous mapping into the real line.
\newblock {\em Topology Appl.}, 12:181--185, 1981.

\bibitem{Todorcevic:1984}
S.~{Todor\v cevi\'c}.
\newblock Trees and linearly ordered sets.
\newblock In {\em Handbook of {S}et-{T}heoretic {T}opology (Kenneth Kunen and
  Jerry~E. Vaughan, eds.)}, pages 235--293. North-Holland, Amsterdam, 1984.

\end{thebibliography}
\end{document}